\documentclass[12pt]{article}

\usepackage{amsmath,amssymb,amsfonts, amsthm,epsf,epsfig}
\usepackage{euscript,graphicx}
\usepackage{longtable}
\usepackage{multirow}

\usepackage{float}
\usepackage{comment}

\usepackage{hyperref}

\newcommand{\ZZ}{\mathbb{Z}}

\theoremstyle{definition}
\newtheorem{theorem}{Theorem}[section]
\newtheorem{lemma}[theorem]{Lemma}

\newtheorem{proposition}[theorem]{Proposition}
\newtheorem{corollary}[theorem]{Corollary}

\newtheorem{remark}[theorem]{Remark}

\newcommand{\dub}{\succ}
\newcommand{\es}{\emptyset}

\newcommand{\Qed}{\hbox{ }\hfill\rule{2.5mm}{3mm}\medskip}

\newcommand{\M}{{\mathcal{M}}}
\newcommand{\N}{{\mathcal{N}}}
\newcommand{\T}{{\mathcal{T}}}
\newcommand{\G}{{\mathcal{G}}}
\newcommand{\cP}{\mathcal{P}}

\newcommand{\cra}{{\leftrightarrow}}

\begin{document}
\begin{center}
{\bf\large
Flag Bicolorings, Pseudo-Orientations, and Double Covers of Maps.
}
\medskip

Hiroki Koike\\
\texttt{yanorado@gmail.com} \\
Institute Andrej Maru\v{s}i\v{c}\\
University of Primorska\\
Koper, Slovenia\\
\medskip

Daniel Pellicer\\
\texttt{pellicer@matmor.unam.mx} \\
Centro de Ciencias Matem\'aticas\\
Universidad Nacional Aut\'onoma de M\'exico\\
Morelia, Mexico \\
\medskip

Miguel Raggi\\
\texttt{mraggi@gmail.com} \\
Centro de Ciencias Matem\'aticas\\
Universidad Nacional Aut\'onoma de M\'exico\\
Morelia, Mexico \\
\medskip

Steve Wilson\\
\texttt{stephen.wilson@nau.edu}\\
Department of Mathematics and Statistics\\
Northern Arizona University\\
Flagstaff, USA \\

\end{center}

\begin{abstract}
This paper discusses consistent flag bicolorings of maps and maniplexes, in their own right and as generalizations of orientations and pseudo-orientations.  Furthermore, a related doubling concept is introduced, and relationships between these ideas are explored.
\end{abstract}

\section{Introduction}

The main goal of this paper is to develop a general theory of flag bicolorings and the related concepts of coverings and pseudo-orientations.

The idea of consistent colorings of the flags of a map with two colors has appeared previously in literature in different contexts. For instance, when the automorphism group of a map which is a polytope has two orbits on the flags, we may color the flags with two colors in such a way that flags in different orbits have different colors. These are called 2-orbit maps in \cite{H2}, but this investigation requires the extra property that the maps be abstract. Specific $k$-colorings of flags are equivalent to the concept of $\mathcal{T}$-{\em compatible} maps introduced in \cite{OPW}, where $\mathcal{T}$ is a class of $k$-orbit maps. This concept, generalized to hypermaps is called $\mathcal{T}$-{\em conservative} in  \cite{B}.  In all these instances, the colorings are used as a tool to work with automorphisms of maps.

Although some of the motivations behind the bicoloring of flags come from  maps admitting symmetries, the ideas can be applied to general maps. In this work we shall not make assumptions on the automorphism groups of the maps in consideration.

Maps admitting a bicoloring of flags in which adjacent flags have different colors are precisely maps on orientable surfaces (see Proposition \ref{oble012} ).  Although the other types of bicolorings of flags that we consider in this paper do not have a known topological equivalence, some results on orientability do translate directly to analogous results on bicolorability. In this sense we may think of bicolorings of flags as a generalization of orientability of maps.

The paper is organized as follows. Sections \ref{sec:maps}, \ref{sec:Ori}, \ref{sec:col} introduce the concepts of maps, pseudo-orientations and bicolorings, respectively. A relationship between pseudo-orientability and bicolorings is provided in Section \ref{sec:ColOri}. Section \ref{sec:op} explores the impact on bicolorability of operations on maps. In Section \ref{sec:doubles}, we discuss a natural double cover of every map which is not bicolorable. The set of bicolorings of any map can be given a natural group structure. In Section \ref{sec:everygroup} we determine, for each possible group, which surfaces admit maps with the given group of bicolorings. Finally, in Section \ref{sec:mpx}, we generalize some of the results about maps to  higher dimensional structures.

\section{Maps}\label{sec:maps}
A map $\M$ is, first and foremost, an embedding of a graph (or pseudograph) on a (compact, connected) surface so that the components of the complement of the embedding (called {\em faces}) are topologically open disks.  We can, for example, regard the cube as an embedding of the graph $Q_3$ on the sphere.

The graph $Q_3$ is an example of a graph which is {\em bipartite}, i.e., its vertices can be colored with two colors so that every edge joins vertices of opposite colors. When speaking about maps, we will use the term {\em vertex-bipartite} to describe a map $\M$ whose underlying graph is bipartite. Similarly, we will call $\M$ {\em face-bipartite} provided its faces can be colored with two colors so that each edge separates faces of opposite colors.  We will call $\M$ {\em edge-bipartite} provided that the edges can be colored with two colors so that edges which are consecutive around a face (and hence around a vertex) have different colors.

\begin{figure}[H]
\begin{center}
\includegraphics[height=30mm]{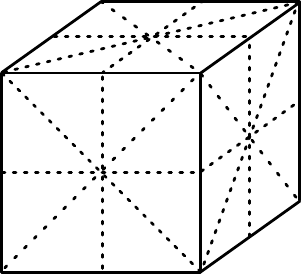}
\caption{The cube divided into flags}\label{fig:cubeflags}
\end{center}
\end{figure}

To look more closely at the structure of a map, we find the following subdivision useful:  choose a point in the interior of each face to call its {\em center} and a point in the relative interior of each edge to be its {\em midpoint}.   Draw dotted lines to connect each face-center with each incidence with the surrounding vertices and edge-midpoints.  The original edges and these dotted lines divide the surface into triangles called {\em flags}. Figure \ref{fig:cubeflags} shows the subdivision of the cube into flags.

\begin{figure}[H]
\begin{center}
\includegraphics[height=30mm]{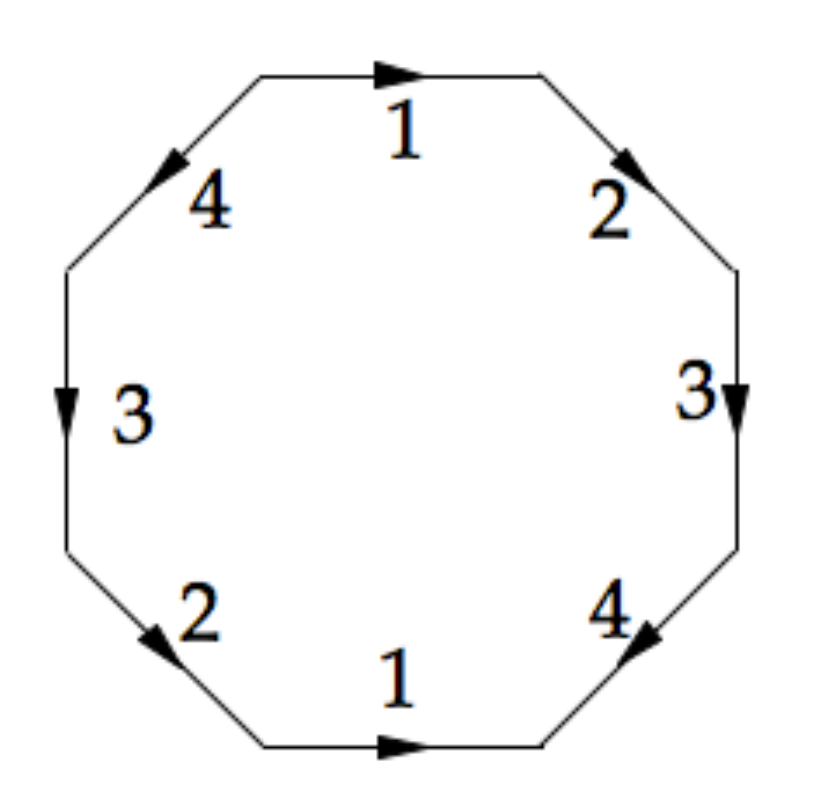}
\hspace{20mm}
\includegraphics[height=30mm]{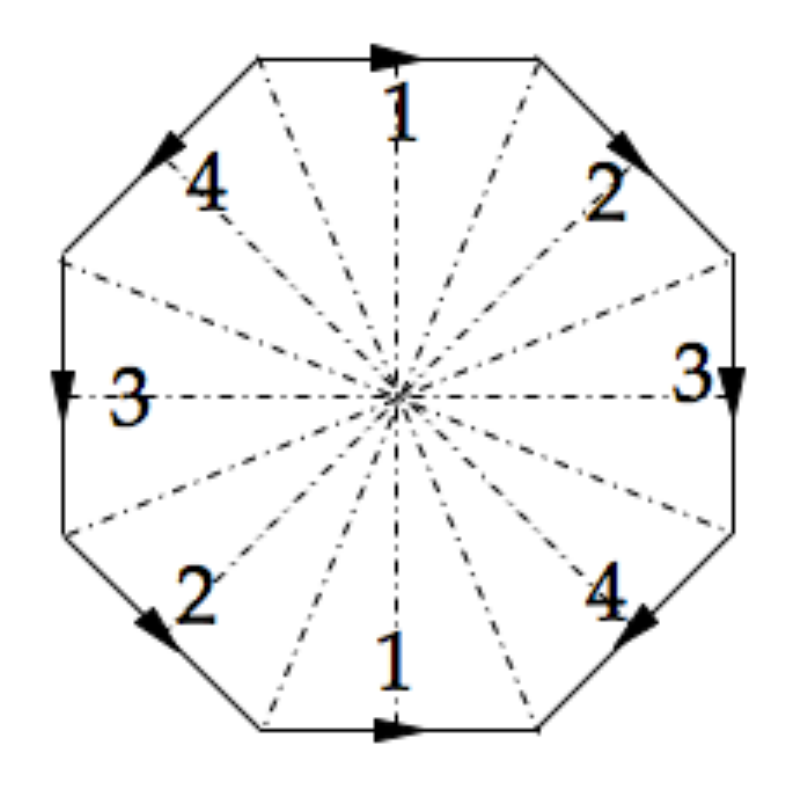}
\caption{The map $M_4$}\label{fig:M4}
\end{center}
\end{figure}

If a face meets a vertex (or an edge) more than once, we emphasize that the face center is connected to each incidence---each appearance---of the vertex or midpoint.  For instance, consider the map $M_4$ shown in the left of Figure \ref{fig:M4}.
This map has only one face, an octagon,  four edges and only one vertex.  Nonetheless, the dissection into flags draws 16 dotted lines, dividing the octagon, and the map, into 16 flags, as shown on the right.

Each flag corresponds to a mutual incidence of face, edge, and vertex, though different flags may correspond to the same triple.

\begin{figure}[H]
\begin{center}
\includegraphics[height=40mm]{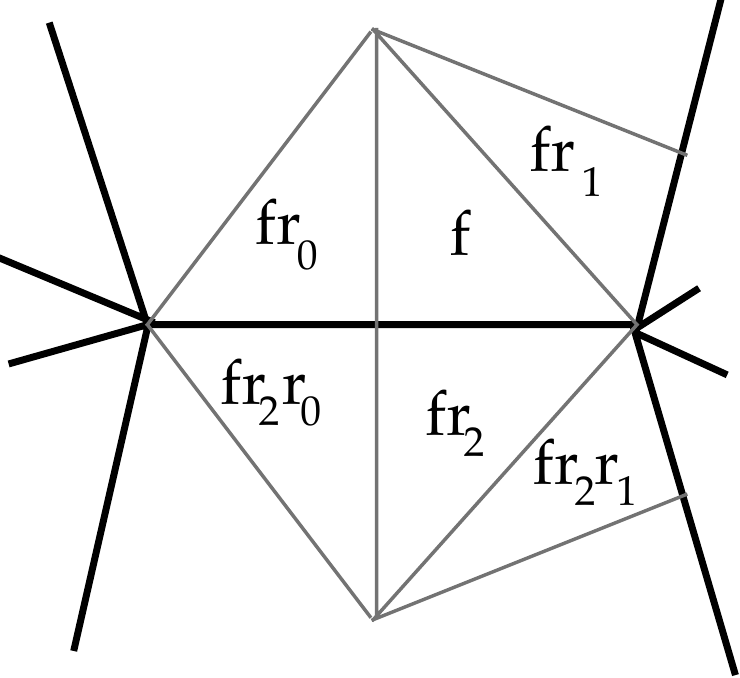}
\caption{Flags in  a map}\label{fig:Flags}
\end{center}
\end{figure}

Let $\Omega$ be the set of flags.  Then let $r_0, r_1, r_2$ be the permutations on $\Omega$ which match each flag $f$ with its three immediate neighbors, as in Figure \ref{fig:Flags}.

In that figure, we see that $f$ and $fr_0$ are adjacent along a face-center-to-edge-midpoint line.  Thus $f$ and $ fr_0$ differ only in their incidences to a vertex, a 0-dimensional face of $\M$.  Similarly, $f$ and $ fr_1$ differ only in their incidences to an edge, a 1-dimensional face, while $f$ and $ fr_2$ differ only in their incidences to a 2-dimensional face. Notice from Figure \ref{fig:Flags} that the flag $r_2$-adjacent to $fr_0$ is also $r_0$-adjacent to $fr_2$.  In other words, as permutations on $\Omega, r_0$ and $r_2$ commute.

We can take a slightly more abstract point of view by defining a map to be a pair $(\Omega, [r_0, r_1, r_2])$ where $\Omega$ is a set of things called {\em flags}, the $r_i$'s and $r_0 r_2$ are fixed point free permutations of order 2 on $\Omega$, the {\em connection group} $C(\M) =\left<r_0, r_1, r_2\right>$ is transitive on $\Omega$, and $r_0$ and $r_2$ commute.  This $C(\M)$  is often called the {\em monodromy group} of the map.   We can then think of vertices in $\M$ as orbits of $\left<r_1, r_2\right>$ in $C(\M)$.   Similarly, edges correspond to orbits of $\left<r_0, r_2\right>$   and faces to orbits of $\left<r_0, r_1\right>$.

We will use the word {\em kite} or {\em corner} for the area within a face where two consecutive edges meet.  More formally, a kite is the union of two flags which are $r_1$-adjacent.

If $\M$ and $\N$ are maps on surfaces $S_{\M}$ and $S_{\N}$, a {\em projection} from $\N$ to $\M$ is a function $\phi$ mapping $S_{\N}$ to $S_{\M}$ which is locally a homeomorphism at all points of $\N$ except perhaps at vertices and/or face-centers, and which sends faces to faces, edges to edges, and vertices to vertices.

In combinatorial terms, if $\M = (\Omega, [r_0, r_1, r_2])$ and $\N = (\Omega', [s_0, s_1, s_2])$, a projection of $\N$ to $\M$ is a function $\phi$ mapping $\Omega'$ to $\Omega$ such that $s_i\phi = \phi r_i$ for all $i \in N$.

We call such an $\N$ a {\em cover} of $\M$.  Notice that if $\N$ is a cover of $\M$, and the pre-image of some one flag in $\M$ has size $k$, then the projection $\phi$ is $k$-to-1 onto every flag of $\M$.  We say then that $\N$ is a {\em $k$-fold} cover of $\M$.

\section{Orientations and Pseudo-Orientations} \label{sec:Ori}
\begin{figure}[H]
\begin{center}
\includegraphics[height=30mm]{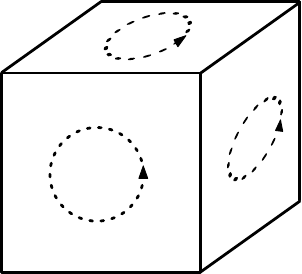}
\caption{Face orientation of the cube}\label{fig:FaceOri}
\end{center}
\end{figure}

We call a map {\em orientable} provided that the surface on which it is embedded is itself orientable.  We check that a map is orientable by giving it a {\em face orientation}; this is an assignment of a circular arrow to each face of $\M$ such that at every edge, the arrows on the faces joined by the edge point along the edge in {\em opposite} directions, as in Figure \ref{fig:FaceOri}.  We can define a {\em vertex orientation} similarly, and it is clear that $\M$ has a vertex-orientation if and only if it has a face-orientation, and this happens if and only if $\M$ is orientable.

\begin{figure}[H]
\begin{center}
\includegraphics[height=42mm]{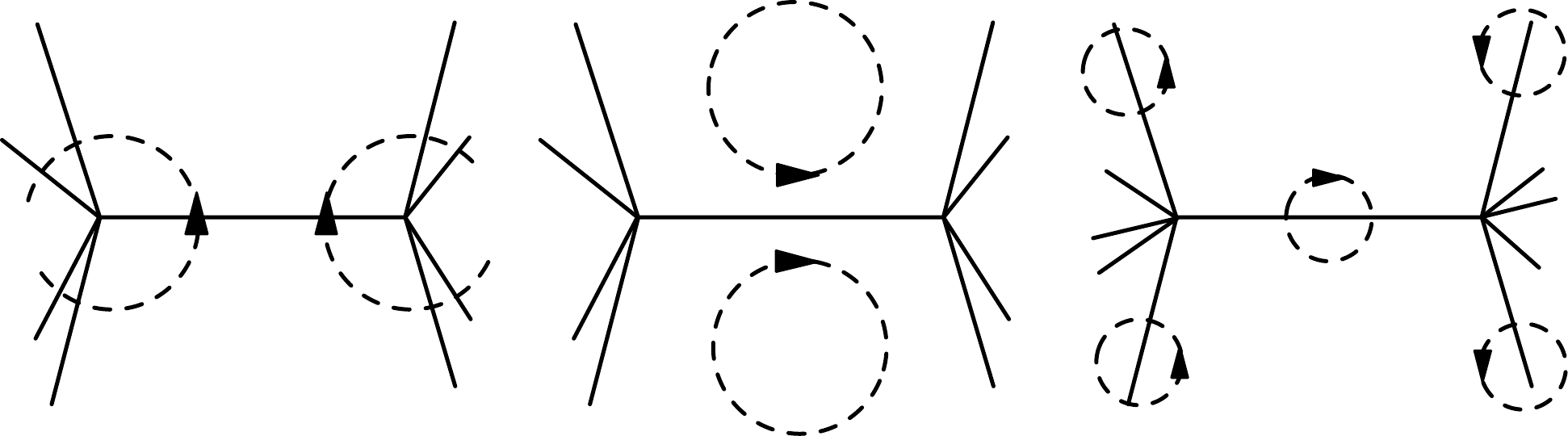}
\caption{Vertex, Face and Edge Pseudo-orientations} \label{fig:PSO}
\end{center}
\end{figure}

The word {\em pseudo-orientation}  has been used in two different ways.  We will use the term {\em vertex pseudo-orientation} (VPSO for short) for what is called in  \cite{WRie} simply a pseudo-orientation.  Here, we mean an assignment of one circular arrow to each vertex so that at each edge the two arrows cross in the {\em  same} direction, as in Figure \ref{fig:PSO}. If we visualize a gear wheel at each vertex so that cogs on the wheels of adjacent vertices mesh, the map is VPSO if we can turn  one wheel, causing  all wheels to turn at the same time. 

In \cite{WRie}, the idea is used to make an important distinction about  $k$-fold  rotary covers of a rotary map $\M$ for which the branching is totally ramified  at vertices.  If $\M$ is non-orientable, $k$ can be larger than 2 only if $\M$ is vertex pseudo-orientable.

Similarly, a {\em face pseudo-orientation} (or FPSO) is an assignment of one circular arrow to each face, as in Figure \ref{fig:PSO}, so that at each edge, the arrows in adjacent faces flow along the edge in the {\em same} direction.  We can simplify this by orienting each edge so that the cycle of edges around each face is consistently oriented.

In \cite{HC} this is called a pseudo-orientation, and the result in that paper is that the {\em Dart Graph} of $\M$ is connected if and only if $\M$ is not face pseudo-orientable.

Finally, we may define an \emph{edge pseudo-orientation} (or EPSO) to be an orientation of the edges such that in every face the direction of the arrows in adjacent edges flows in the \emph{same} direction (either into or out of the face) as in Figure \ref{fig:PSO}.

\section{Colorings} \label{sec:col}
\begin{figure}[H]
\begin{center}
\includegraphics[height=35mm]{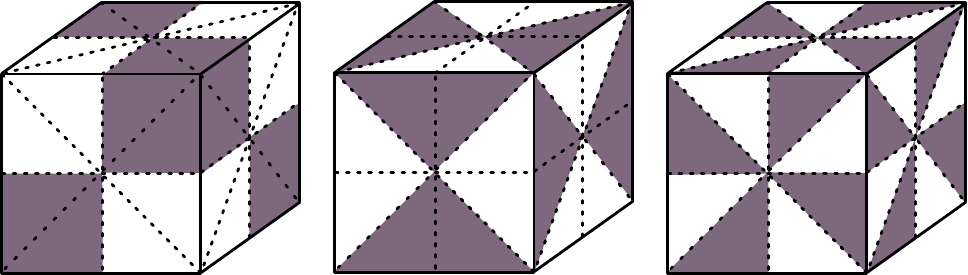}
\caption{Colorings of the cube} \label{fig:CubeCols}
\end{center}
\end{figure}

In this paper, the word `coloring' will be used to describe what might be more fully notated as a `consistent flag 2-coloring'.  If $I$ is any subset of $R = \{0, 1, 2\}$, then an $I$-coloring of a map $\M$ is a function $a: \Omega \rightarrow \ZZ_2$ such that for every flag $f$, if $j \in I$, then $a(fr_j) \neq a(f)$, while  if $j \notin I$, then $a(fr_j) = a(f)$. To say that in another way,

$$
a(fr_j)+ a(f)  = \left\{ \begin{array}{ll}
         1 & \text{if } j \in I\\
         0 & \text{if } j \notin I.\\
    \end{array} \right.
$$
Figure \ref{fig:CubeCols} shows a $\{0\}$-coloring, a $\{1,2\}$-coloring, and an $R$-coloring of the cube. The reader should check that the colorings extend unambiguously to the unseen faces of the cube as well.

Let us first notice that if $a$ is an $I$-coloring of $\M$, then so is $1-a$, and these are, in fact, the only $I$-colorings of the map.  

Give the name {\em sub-edge} to the segments forming sides of the flags.  These are of three types:

\begin{description}
  \item[Type 0]: midpoint-of-edge to face-center,
  \item[Type 1]: face-center to vertex,
  \item[Type 2]: vertex to midpoint-of-edge.
\end{description}

Each flag $f$ is adjacent to flag $fr_i$ along their mutual sub-edges of type $i$. If $\M$ is $I$-colorable, therefore, the sub-edges bounding each color are the sub-edges of types in $I$.

Note that when considering maps in class $2_I$ as defined in \cite{H2}, the subset $I$ is used to indicate that $i$-adjacent flags belong to the {\em same} orbit. This means that maps in class $2_I$ admit an $(R \setminus I)$-coloring. In the context of colorings it turns out to be more practical to define $I$-coloring requiring $i$-adjacent flags to be a different color rather than requiring $i$-adjacent flags to be the same color if and only if $i \in I$.

\subsection{The group \texorpdfstring{$T(\M)$}{T(H)}}

Let $\Delta$ stand for the symmetric difference operation:  $I\Delta J =( I\cup J)\backslash(I\cap J)$. Let $\cP$ be the power set of $R$. It is well-known that $\cP$ is a group under $\Delta$ and it is an elementary Abelian group of order eight whose identity element is $\emptyset$.

\begin{proposition}\label{Delta}
If a map $\M$ admits an $I$-coloring and a $J$-coloring, then it admits an $I\Delta J$-coloring.
\end{proposition}
Proof:  It is easy to verify that the sum of an $I$-coloring and a $J$-coloring is an $I\Delta J$-coloring \Qed

From this proposition we see that, for each map $\M$, the sets $I$ for which $\M$ has an $I$-coloring form a set $T(\M)$ which is a subgroup of $\cP$.

\subsection{Flag-colorings and map-colorings}
Whenever $I$ is a subset of size 1, $I$-colorable maps have a nice characterization.
\begin{proposition}
\label{Vbip}
A map $\M$ is vertex-bipartite if and only if $\M$ has a $\{0\}$-coloring.
\end{proposition}
Proof: Given such a coloring, all the flags incident with one vertex will be the same color, while those at an adjacent vertex will be the other color.  This is a bipartite coloring of the vertices. Conversely, given a bipartition of the vertices, color all flags incident with black vertices with color 1 and all those incident with white vertices with color 0. This is then a  $\{0\}$-coloring of $\M$.\Qed

Similar proofs lead us to these:
\begin{proposition}
\label{Fbip}
A map $\M$ is face-bipartite if and only if $\M$ has a $\{2\}$-coloring. (\textit{i.e.} a $\{2\}$-coloring corresponds to a bipartite coloring of the faces of $\M$.)
\end{proposition}

\begin{proposition}
\label{Ebip}
A map $\M$ is edge-bipartite if and only if $\M$ has a $\{1\}$-coloring. (\textit{i.e.} a $\{1\}$-coloring corresponds to a bipartite coloring of the edges of $\M$.)
\end{proposition}

\subsection{Colorings and words}\label{words}
Consider a {\em cycle} $f, fr_{i_1}, fr_{i_1}r_{i_2},fr_{i_1}r_{i_2}r_{i_3}, \dots,fr_{i_1}r_{i_2}r_{i_3}\dots r_{i_k} = f$ of flags in a map $\M$.   We abbreviate that $(f, W)$, where $W$ is the string $i_1i_2i_3 \dots i_k$.  Write $r_W
$ for $r_{i_1}r_{i_2}r_{i_3}\dots r_{i_k}$. We call a cycle $(f, W)$ \emph{ $I$-consistent} if the number of occurrences  in $W$ of indices which are in $I$ is even.  If $(f, W)$ is $I$-consistent, this says that, at least along the cycle, we can color flags from two colors so that flags which are $i$-connected for $i \notin I$ are the same color and those for which $i \in I$ are not. A cycle which is not $I$-consistent is \emph{$I$-inconsistent}. It follows, then that $\M$ is $I$-colorable if and only if every cycle in $\M$ is $I$-consistent.

Consider this quite general fact: for any function $F$  from a set $X$ to $\ZZ_2$, extend $F$ to the power set of $X$ by $F(A) = \Sigma_{x \in A}F(x)$ for each subset $A$ of $X$. Then for any subsets $A, B$ of $X$, it is clear that $F(A\Delta B)= F(A)+F(B)$.

\begin{lemma} \label{lemma13}
For any subsets $I, J$ of $R$, if $K = I\Delta J$, then every cycle in $M$ is consistent either for exactly one of $I, J, K$ or for all three.
\end{lemma}

Proof: Fix a cycle $(f, W)$ of $\M$, and define $F$ on $R$ with $F(i)$ being the parity of the number of occurrences in $W$ of $r_i$; furthermore, consider $F$ extended to $\cP$.  The cycle is $I$-consistent, then, provided that $F(I) = 0$ .  Now consider $F(I) + F(J)=F(I\Delta J) = F(K)$. These are elements of $\ZZ_2$, and so we can rephrase that in this form: $F(I)+F(J)+F(K) = 0$.  The number of zeros in $\{F(I),F(J),F(K)\}$ must be 1 or 3;  the conclusion follows directly.   \Qed

We will have use for this lemma in Section \ref{sec:doubles}.

\subsection{Connected sum of maps}\label{sec:connectedsum}

A natural operator in topology is the \emph{connected sum} of two surfaces, which is defined to be the surface formed by cutting a small disk from each of the two surfaces and attaching them along the newly created border. We extend this definition to maps by causing the attachment sets to be the boundaries of two suitably chosen faces.

Formally, let $\M = (\Omega, [r_0, r_1, r_2])$ and $\N = (\Omega', [s_0, s_1, s_2])$ be maps and suppose $f_\M$ and $f_\N$ are flags of these maps for which their respective faces $F_\M$ and $F_\N$ have the same number $k$ of sides. Furthermore, suppose that no flag of $F_\M$ is $r_2$-connected to any of the flags of  $F_\M$, and similarly for $F_\N$. Then we define the \emph{connected sum} $\M\oplus \N$ with respect to $(f_\M,f_\N)$ to be the map whose flag set is   $(\Omega\backslash F_\M)\cup(\Omega'\backslash F_\N)$, with connections $t_0, t_1, t_2$, where $t_0$ and $t_1$ are the restrictions of  $r_0\cup s_0,  r_1\cup s_1$to this set, and $t_2$ is the same except that, for each $j, f_\M(r_0r_1)^jr_2$ is $t_2$-connected to  $f_\N(s_0s_1)^js_2$ and $ f_\M r_1(r_0r_1)^jr_2$ is $t_2$-connected to  $f_\N s_1(s_0s_1)^js_2$.  Because of the restriction on the faces, these flags are both in the new flag set, and it is easy to check that $t_0t_2 = t_2t_0$.

The following Figure illustrates this construction:

\begin{figure}[H]
\begin{center}
\includegraphics[scale=1]{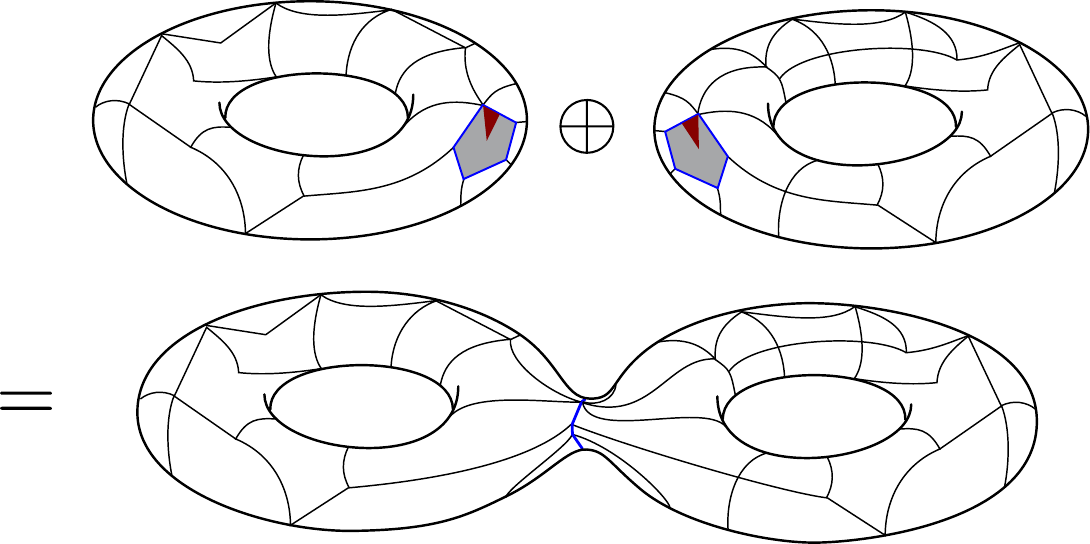}
\caption{A connected sum of two maps}
\label{fig:connectedsum}
\end{center}
\end{figure}

\begin{proposition}
\label{th:oplus}
  Suppose $\M$ and $\N$ are $I$-colorable maps for some $I\subseteq R$ and suppose $\M\oplus \N$ exists for some flags $f_\M$ and $f_\N$. Then $\M\oplus \N$ is $I$-colorable.
\end{proposition}
Proof: Choose an $I$-coloring for $\M$ and one for $\N$ in such a way that $f_\M$ and $f_\N$ have opposite colors if $2\in I$ and the same color if $2 \notin I$. Then the induced coloring provides the desired coloring on $\M\oplus\N$. \Qed

\vspace{0.3cm}

As we'll see in the following theorem, it is almost always the case that the group $T(\M \oplus \N)$ is equal to $T(\M) \cap T(\N)$. This can fail, roughly, when for some $I \subseteq R$, $\M\setminus F_\M$ is $I$-colorable but $\M$ is not, or similarly for $\N$.

\begin{theorem}\label{t:plus}
    Suppose $\M$ and $\N$ are maps and suppose $\M\oplus \N$ exists for some flags $f_\M$ and $f_\N$. Furthermore, suppose that for each $I\notin T(\M)$ there exists an $I$-inconsistent cycle that doesn't include any of the flags of face $F_\M$, and similarly for $\N$. Then $T(\M\oplus \N) = T(\M) \cap T(\N)$.
\end{theorem}

Proof:  For  $I \subseteq R$, if $I\in T(\M) \cap T(\N)$, then $I\in T(\M\oplus \N)$ by Proposition \ref{th:oplus}.  If $I\notin T(\M) \cap T(\N)$, then without loss of generality, assume that $\M$ is not $I$-colorable.  Then $\M$ has cycles which are $I$-inconsistent, and by the hypothesis, at least one of them includes none of the flags in $F_\M$.  This cycle, then, occurs in $\M\oplus \N$, and is still $I$-inconsistent.  Thus $I \notin T(\M\oplus \N)$.  \Qed

\section{Colorings and Orientations}
\label{sec:ColOri}
We now offer and prove a series of propositions connecting the ideas of orientations and colorings. First, we show a collection of easy facts to embody our belief that colorings generalize orientations. We then generalize these facts in Theorem \ref{thm:super}.

\begin{proposition}\label{oble012}
A map is orientable if and only if it has an $R$-coloring.
\end{proposition}
Proof:  If $\M$ has an $R$-coloring, orient each face so that, along each edge, the orientation points from the flag with color 0 to the flag with color 1.  This is consistent within the face and faces that meet along an edge have orientations which meet correctly.  Conversely, given the orientation, assign colors so that along each edge within a face, the orientation points from the flag with color 0 to the flag with color 1.  This causes colors to alternate within each face.  Because the orientation opposes the orientation in each adjacent face, the colors alternate at each edge, giving an $R$-coloring.\Qed

Similar arguments give us the following two propositions:

\begin{proposition}\label{th:fpso01}
 A map $\M$ has a face pseudo-orientation if and only if $\M$ admits a $\{0,1\}$-coloring.
\end{proposition}

\begin{proposition}\label{vpso12}
A map $\M$ has a vertex pseudo-orientation if and only if $\M$ admits a $\{1,2\}$-coloring.
\end{proposition}

\subsection{\texorpdfstring{$I$}{I}-Pseudo-Orientability}\label{sec:IPSO}

Our aim in this subsection is to generalize the definition of pseudo-orientability and the previous results. For any subset $I$ of $R$ (except $\emptyset$), we form the map $X = X(\M, I)$ whose faces are the regions formed by conjoining flags which are connected by each $r_i$ for which $i \notin I$. For example, if $\M$ is the cube, then $X(\M,\{1,2\})$ and $X(\M,\{0,2\})$ are shown below.

\begin{figure}[H]
\begin{center}
\includegraphics[scale=1]{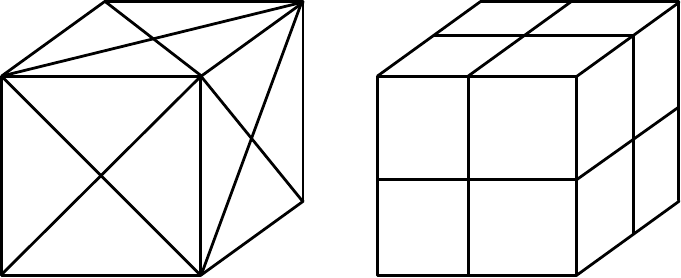}
\caption{$X(\M,\{1,2\})$ and $X(\M,\{0,2\})$}
\end{center}
\end{figure}

Notice that the boundaries of the faces in $X(\M, I)$ are the sub-edges of types that are {\em in} $I$.  Thus $\M$ is $I$-colorable if and only if  $X(\M, I)$ is face-bipartite. However, the map $X(\M, I)$ can be constructed even if $\M$ is not $I$-colorable.

We call $\M$ \emph{$I$-pseudo-orientable} provided that $X(\M,I)$ is face pseudo-orientable. If this happens, we can orient the edges of $X$ so that around each face of $X$, the arrows all point the same way. We say a map is $\emptyset$-pseudo-orientable when it is orientable.

It is clear that $X(\M,\{2\})=\M$ and so $\{2\}$-PSO is equivalent to FPSO. Observation \ref{th:fpso01} shows that this happens exactly when $\M$ is $\{0,1\}$-colorable. Similarly, $\M$ is $\{0\}$-PSO when it is VPSO, and this happens exactly when $\M$ is $\{1,2\}$-colorable. Surprisingly, this generalizes completely:

\begin{theorem}\label{thm:super}
  Let $\M$ be a map. Then $\M$ is $I$-pseudo-orientable if and only if $\M$ is $(R\backslash I)$-colorable.
\end{theorem}

Proof:   First, recall that the sub-edges of $\M$ are of three different types:
\begin{description}
  \item[Type 0]: midpoint-of-edge to face-center,
  \item[Type 1]: face-center to vertex,
  \item[Type 2]: vertex to midpoint-of-edge.
\end{description}

Furthermore, let us think of the sub-edges as directed. Call the directions above  the `standard' orientations for the sub-edges, and the reverse of these are the non-standard orientations.

In order to form $X(\M,I)$, we delete or ignore sub-edges of each type $i\notin I$, so we only consider sub-edges of type in $I$. Suppose $\M$ is already $(R\backslash I)$-colored. Then in $X(\M, I)$, both flags around each sub-edge have the same color. Then we may assign this as the color of the sub-edge.  Now give every sub-edge with color 0 the standard orientation, and each sub-edge with color 1 the non-standard direction.

If $I$ has 2 elements, then consecutive sub-edges around a face of $X$ of different types are the same color, while those of the same type are opposite colors.  This causes the directions to be consistent about a face.  On the other hand, if $I$ has only one element, then the sub-edges around a face of $X$ alternate in color and are all of the same type; again, this forces the orientation to be consistent.

The following examples show how this works:

\begin{enumerate}
  \item If $I=\emptyset$, then by definition the map is $I$-pseudo-orientable if and only if it is orientable, and this happens if and only if it is $R$-colorable.
  \item If $I=\{2\}$, then $X(\M,I)=\M$, which is face pseudo-orientable if and only if it is $\{0,1\}$-colorable. In this case, we orient black type 2 sub-edges from vertex to midpoint-of-edge but white type 2 sub-edges from midpoint-of-edge to vertex.
  \item If $I=\{1,2\}$, we orient black type 1 sub-edges from vertex to face-center and black type 2 sub-edges from midpoint-of-edge to vertex. We orient the opposite colors in the opposite direction:
    \begin{center}
      \includegraphics[scale=1]{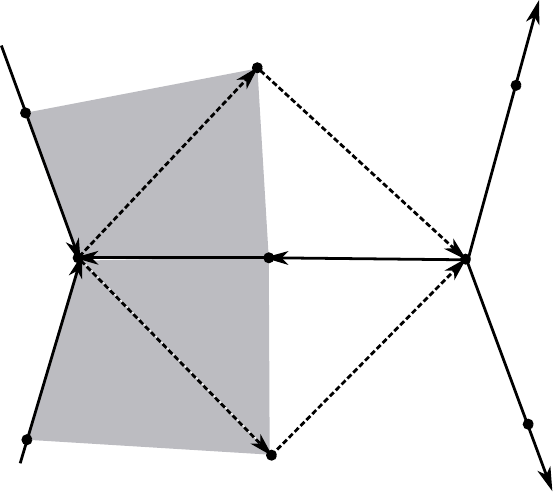}
    \end{center}
  \item If $I=R$, the  $\emptyset$-coloration colors all sub-edges the same and so the faces of $X$, the flags of $\M$, are all oriented consistently.  Every map is $\emptyset$-colorable and every map is $R$-PSO.
  \begin{center}
      \includegraphics[scale=1]{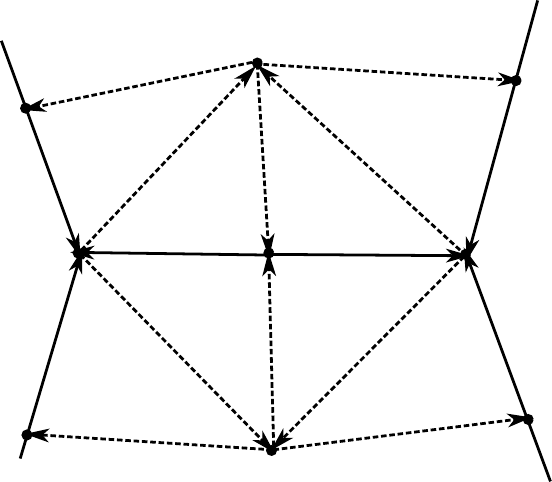}
  \end{center}
 
\end{enumerate}

For the converse we may use the orientation to color in the way prescribed to match the previous list. \Qed

\subsection{The Cheat Sheet}

We summarize here the results of previous sections relating to the question of which maps have $I$-colorings for a given $I$.

\begin{center}
\begin{tabular}{|r|l|c|}
\hline
$I$&$\M$&Requires:\\
\hline\hline
$\emptyset$& All flags have the same color & -\\
\hline
$\{0\}$ & Vertex-bipartite & Faces have even degree\\
\hline
$\{1\}$ & Edge-bipartite & Faces and vertices have even degree\\
\hline
$\{2\}$ & Face-bipartite & Vertices have even degree\\
\hline
$\{0,1\}$ & Face Pseudo-Orientable & Vertices have even degree\\
\hline
$\{0,2\}$ & Edge Pseudo-Orientable & Faces and vertices have even degree\\
\hline
$\{1,2\}$ & Vertex Pseudo-Orientable & Faces have even degree\\
\hline
$R$ & Orientable & - \\
\hline
\end{tabular}
\end{center}

\section{Operators}
\label{sec:op}

In this section we shall see how colorings and orientations interact with common map operators such as dual, Petrie, opposite and medial.

Let $\M = (\Omega, [r_0, r_1, r_2])$ be any map.  Then $D(\M)$ (or the \emph{dual} of $\M$) is the map $(\Omega, [d_0, d_1, d_2])$, where $d_j = r_{2-j}$.  This corresponds exactly to the classical geometric dual of a polyhedron. The maps $\M$ and $D(\M)$ lie on the same surface.

The {\em Petrie} of $\M$, $P(\M)$, is the map $(\Omega, [p_0, p_1, p_2])$, where $p_0 = r_0r_2, p_1 = r_1$, and $p_2 = r_2$.  This is a less familiar operator on maps. The faces of $P(\M)$ are the Petrie paths (see \cite{Wop} ) of $\M$ and vice-versa. Note that the vertices and the edges are preserved by the operation.

We call a map formed from $\M$ by any composition of the operation $D$ and $P$ a {\em  direct derivate} of $\M$. Of special interest among the direct derivates of $\M$ is $opp(\M) = PDP(\M) = DPD(\M$).  Formally, this is $(\Omega, [s_0, s_1, s_2])$ where $s_0 = r_0, s_1 = r_1$ and $s_2 = r_0r_2$.  More intuitively, $opp(\M)$  is formed from $\M$ by cutting along the edges and then re-attaching along the same matching edges but  with the reverse local orientation.  See \cite{Wop} for more information about these operators.

The {\em medial} of a map $\M$ is drawn on the same surface as $\M$. The vertices of the medial are the edge-midpoints of $\M$ and two are connected  by an edge diagonally across each kite to which both belong.

\begin{figure}[H]

\begin{center}
\includegraphics[height=40mm]{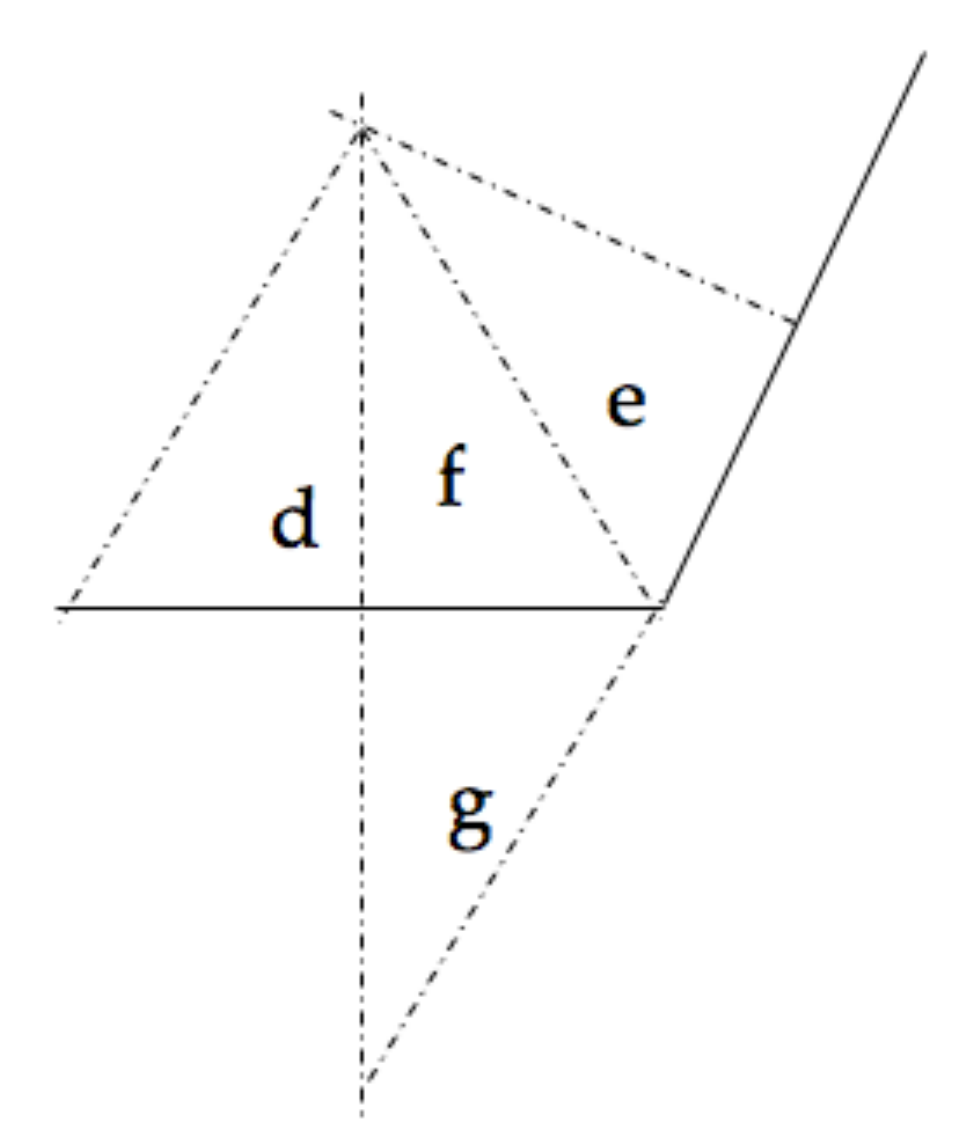}
\includegraphics[height=40mm]{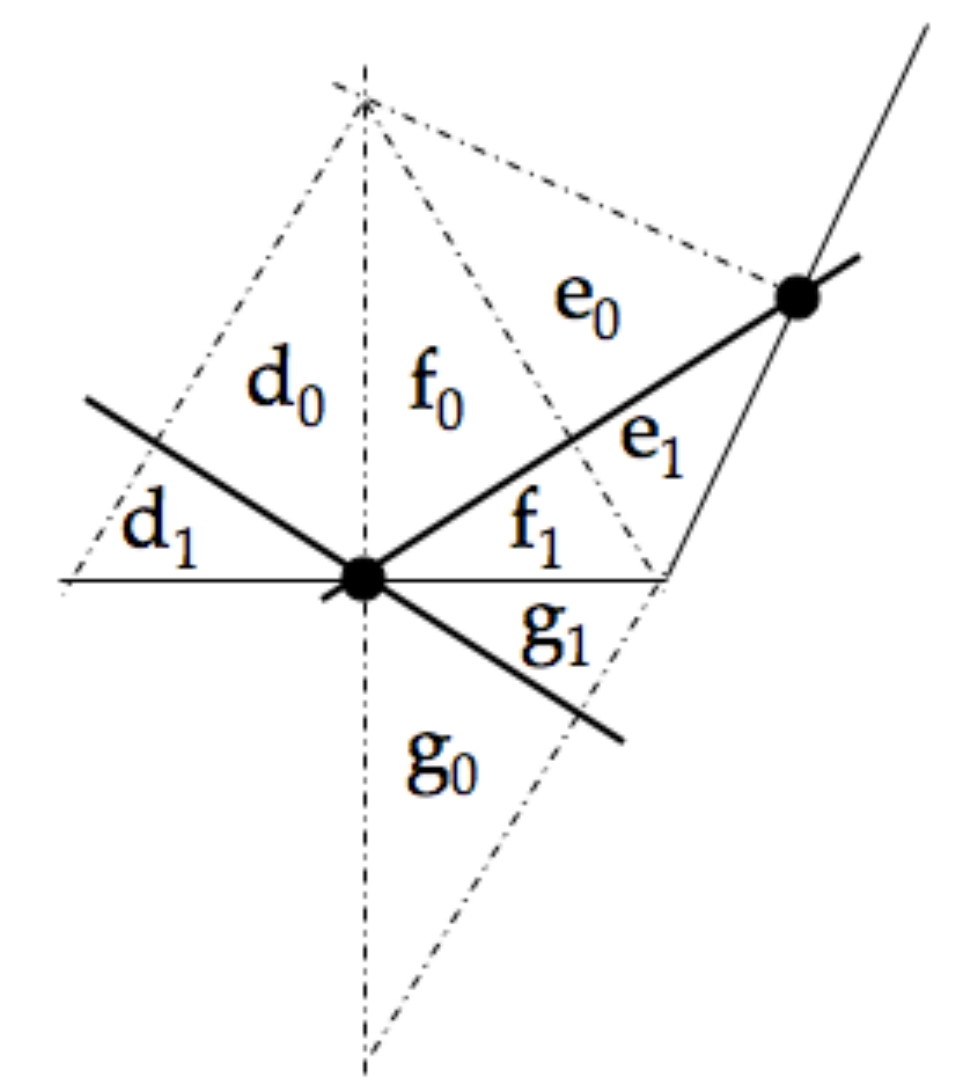}
\caption{Flags in a map and its medial}
\label{fig:Medial}
\end{center}
\end{figure}

Figure \ref{fig:Medial} shows, on the left, the flags adjacent to a flag $f$ in a map.  On the right we see the edges of the medial drawn in in heavy lines.  These divide each flag into two.  This suggests the combinatorial form: the medial of $\M  = (\Omega, [r_0, r_1, r_2])$ is the map $(\Omega\times\ZZ_2, [s_0, s_1, s_2])$ where $s_0, s_1, s_2$ are defined by
\begin{enumerate}
\item $(f, i)s_0 = (fr_1,i)$
\item $(f, 0)s_1 = (fr_0, 0)$
\item $(f, 1)s_1 = (fr_2, 1)$
\item $(f, i)s_2 = (f, 1-i)$.
\end{enumerate}

\subsection{Operators and colorings}

In this subsection, we describe the ways in which the colorings of a map come from the colorings of maps related to it by operators.

\begin{lemma}\label{DI}
If $\M$ admits an $I$-coloring, then $D(\M)$ admits a $J$-coloring, where $J = \{2-i: i \in I\}$
\end{lemma}
Proof: This follows from the definition of the dual. \Qed

\begin{proposition}
\label{oppI}
If $\M$ admits an $I$-coloring, then $opp(\M)$ admits a $J$-coloring, where:
  $$J = \left\{ \begin{array}{cc}
    I & \text{ if } 0 \notin I \\
    I\Delta \{2\} & \text{ if } 0 \in I.\\
  \end{array} \right.$$
\end{proposition}
Proof:  Suppose $\M$ is $I$-colored. We shall take opposites but keep the coloring of each flag. Let $f$ be a flag. Then the $r_0$ and $r_1$-adjacent flags in $opp(\M)$ are the same as in $\M$, and the $r_2$ adjacent flag is $fr_2 r_0$. If $0 \notin I$, then $fr_2$ is the same color as $fr_2 r_0$, which means the coloring doesn't change when taking opposites. Analogously, if $0 \in I$, $fr_2$ is the opposite color as $fr_2 r_0$, so taking opposite changes the color of the flag $r_2$-adjacent to $f$.

\begin{proposition}
\label{PI}
If $\M$ admits an $I$-coloring, then $P(\M)$ admits a $J$-coloring, where:
  $$J = \left\{ \begin{array}{cc}
    I & \text{ if } 2 \notin I \\
    I\Delta \{0\} & \text{ if } 2 \in I. \\
  \end{array} \right.$$
\end{proposition}
Proof: This follows from the previous theorem applying Lemma \ref{DI}. \Qed

\begin{corollary}
 A map $\M$ has all colorings, {\em i.e.} $T(\M) = \cP$, if and only if the surfaces of $\M$, $opp(M)$ and $P(M)$ are orientable.
\end{corollary}

Proof: Indeed, if $R \in T(\M)\cap T(P(\M)) \cap T(opp(\M))$, by the preceding theorems we conclude that $\{0,1\}$ and $\{1,2\}$ are both in $T(\M)$, which means $R \Delta \{0,1\}=\{2\}$ and $R \Delta \{1,2\} = \{0\}$ are also in $T(\M)$, and so $T(\M) = \cP$, since $\{\{0\}, \{2\}, R\}$ is a generating set of $\cP$ under $\Delta$.

Now, suppose $T(\M)=\cP$. By previous results, if $\M$ admits every coloring, then $opp(\M)$ and $P(\M)$ admit every coloring.  Thus all three admit an $R$-coloring and so are orientable. \Qed

\begin{theorem}
\label{med}
The medial map of any map can always be $\{2\}$-colored. Furthermore, for each of the pairs shown in the table below, a map $\M$ has an $I$-coloring if and only if the medial of $\M$ has a $J$-coloring.

\begin{center}
\begin{tabular}{|c|c|}
\hline
 $I$ & $J$\\
\hline
$\{1\}$ & $\{0\}$\\
\hline
$\{0,2\}$ & $\{1\}$\\
\hline
$R$ & $R$\\
\hline
\end{tabular}
\end{center}
\end{theorem}

Proof: The medial of any map always admits a $\{2\}$-coloring (face-bipartition), since we have two types of faces: those that surround vertices of $\M$, and those that surround face-centers of $\M$. For the first row, if the map is edge-bipartite, then we may color each vertex of the medial with the color corresponding to the color of the edge it was the midpoint of and vice-versa. For the second, note that an edge of the medial corresponds to a kite of the original. Finally, the medial maps lies on the same surface as the original.\Qed

\section{Doubles}\label{sec:doubles}
Given a map $\M = (\Omega, [r_0, r_1, r_2])$ and a subset $I$ of $R = \{0, 1, 2\}$, we define the {\em $I$-double} of $\M$, denoted $I\dub \M$,  to be the map $\N = (\Omega\times\ZZ_2, [s_0, s_1, s_2])$, where

$$
f_is_j  = \left\{ \begin{array}{ll}
         (fr_j)_{i} & \text{ if } j \notin I\\
          (fr_j)_{1-i} & \text{ if } j \in I.\end{array} \right.
$$
Here, we write ``$f_i$'' for $(f, i)$.

Notice that the function sending $f_i$ to $f$ is a projection of $\N$ onto $\M$, and so $I\dub\M$ is a covering of $\M$.  Also notice that the function which sends  $f_i$ to $i \in\ZZ_2$, is an $I$-coloring of $I\dub\M$.

Now, if $\M$ has an $I$-coloring, then $I\dub \M$ consists of two disjoint copies of $\M$; in this case we discard one copy and say that $I\dub \M$ is isomorphic to $\M$.  If $\M$ does not have an $I$-coloring then  $I\dub \M$ truly is a double cover (i.e., a 2-fold topological covering which may or may not be branched at vertices and/or face-centers) of $\M$.  In this case, for every cycle $(f, W)$ in $\M$, if it is $I$-consistent, then $f_0s_W = f_0$, and $f_1s_W = f_1$, and so the cycle is covered by two cycles of the same length.  And if $(f, W)$ is $I$-inconsistent, then  $f_1s_W = f_0$, and $f_0s_W = f_1$.  Then the cycle $(f, W)$ is covered by a single cycle of twice the length of $(f, W)$.  In that case the covering cycle is $(f_0, W^2)$, which is consistent for any subset of $R$.

We have remarked above that the map $I\dub \M$ always has an $I$-coloring. Moreover,  we claim, it is universal minimal in the sense of the following proposition:

\begin{proposition}
\label{th:mindub}  If $\N$ is any map that is a covering of $\M$ and which is $I$-colorable, then $\N$ must itself cover $I\dub \M$.
\end{proposition}
Proof:  Let $a$ be an $I$-coloring of  $\N$ and let $\phi$ be a projection of $\N$ onto $\M$.  Define $\phi'$ mapping $\N$ onto $I\dub\M$ to send the flag $f$ of $\N$ to the flag $\phi(f)_{a(f)}$ of $I\dub\M$.  Because $a$ is an I-coloring of $\N$, this $\phi'$ is a projection, as required.  \Qed

  An example of a double covering is $\N = R\dub \M$ when $\M$ is a non-orientable map.  Because $\M$ does not have an $R$-coloring, $\N$ is truly a double cover of it, and since $\N$ then {\em does} have an $R$-coloring, it must be orientable and so $\N$ is the orientable twofold cover (some say the {\em canonical} cover) of $\M$.

Another example  of the double construction in the literature comes from Sherk's 1962 paper \cite{Sh}.  At that time few non-toroidal chiral maps were known.  Sherk constructed an infinite family of such maps on an infinite number of different surfaces, using a simple and clever technique.  Starting with a chiral map $\M$ on the torus, with triangles as faces and six of them incident to each vertex, he made a new map $\N$ in the following way.  For each vertex $v$ of $\M$, he created vertices $v_0, v_1$; for each edge $\{u,v\}$, he made edges $\{u_0, v_1\}$ and $\{u_1, v_0\}$; and each triangular face $[u,v,w]$ became a hexagonal face $[u_0, v_1, w_0, u_1, v_0, w_1]$.  This map has hexagonal faces, still meeting 6 at a vertex.  Thus if $\M$ had $D$ vertices, $2D$ faces and $3D$ edges, then $\N$ has  $2D$ vertices, $2D$ faces and $6D$ edges.  It follows that $\N$ must lie on the orientable surface of characteristic $-2D$, i.e., of genus
$D-1$.  From our point of view, we can construct $\N$ as $\{0\}\dub \M$.

\begin{figure}[H]
\begin{center}
\includegraphics[height=45mm]{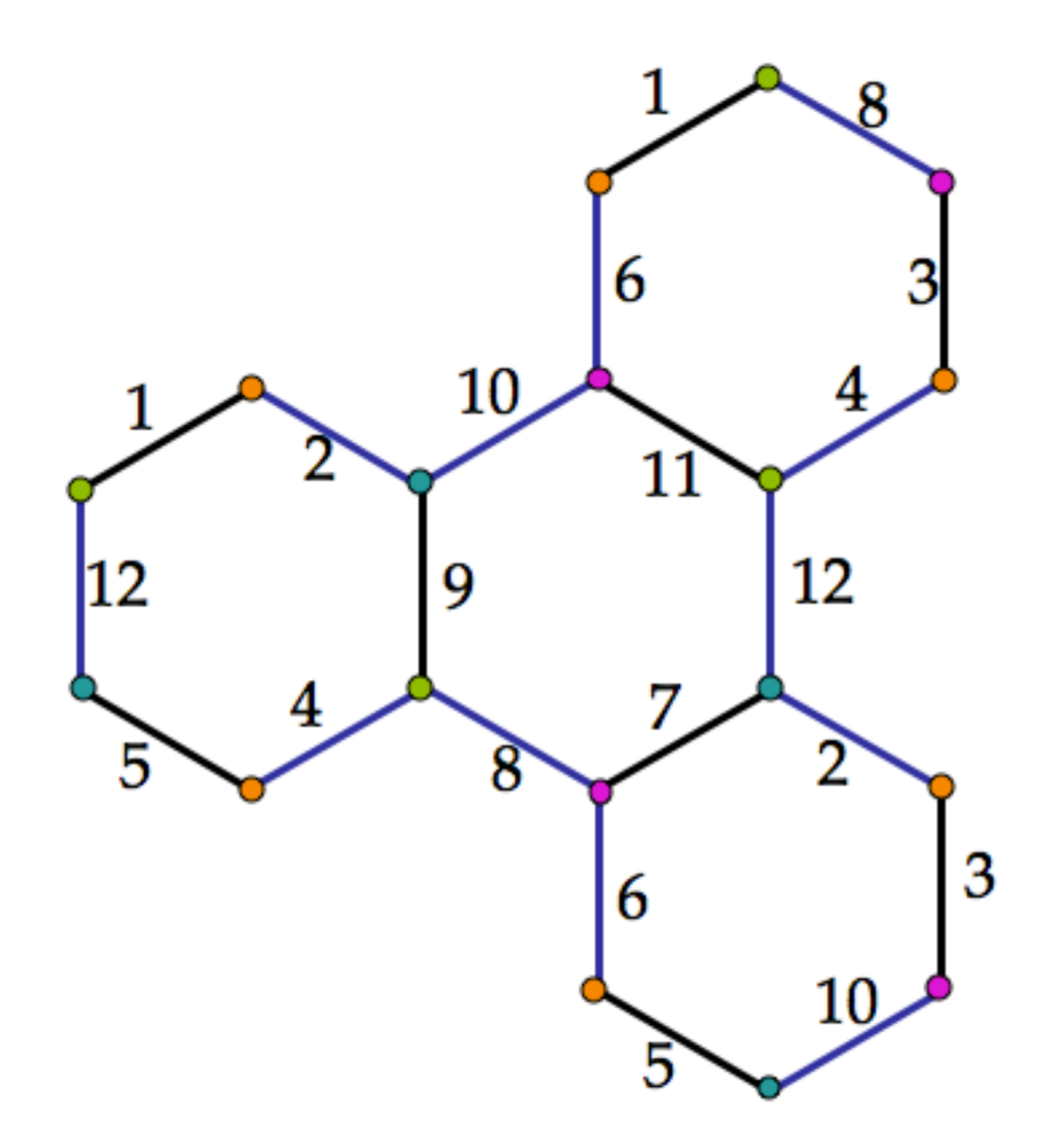}
\caption{The Map $\{1\}\dub\T$}
\label{fig:1YT}
\end{center}
\end{figure}

As a final example, the reader might like to verify that if $\T$ is the tetrahedron then $\{0,2\}\dub\T$ is the orientable map with hexagonal faces, six of them around each vertex, shown in Figure \ref{fig:1YT}.

\subsection{Properties of double covers}
We present here some of the facts about $I$-doubles:

\begin{proposition}
\label{th:Jup}
If $\M$ has a $I$-coloring then so does any $J\dub \M$.
\end{proposition}
Proof:  In fact, by definition {\em any} covering of $\M$ inherits an $I$-coloring from that of $\M$.\Qed

\begin{proposition}\label{t:IDJ}
If $\M$ has a $J$-coloring then any $I\dub \M$ is isomorphic to $(I\Delta J)\dub \M$.
\end{proposition}
Proof:  For any map $\M$, $I\dub \M$ is the smallest covering of $\M$ which admits an $I$-coloring.  If $\M$ has a $J$-coloring, then  $I\dub \M$ has a $J$-coloring as well.  Thus, by Proposition \ref{Delta},  $I\dub \M$ admits an $(I\Delta J)$-coloring, and whether $\M$ has an $I$-coloring or not, this covering of $\M$ must be minimal of those admitting an $(I\Delta J)$-coloring.  Thus $I\dub \M\cong (I\Delta J)\dub \M$. \Qed

The following two results summarize the examples given above.

\begin{proposition}
If $\M$ is non-orientable, then $R \dub \M$ is isomorphic to the orientable double covering of $\M$.
\end{proposition}

We say that a map $\N$ is the {\em Sherk double cover} of a non-vertex-bipartite map $\M$ whenever $\N$ can be obtained from $\M$ by the above procedure used by Sherk to obtain the chiral maps with hexagonal faces.

\begin{proposition}
If $\M$ is not vertex-bipartite, then $\{0\}\dub \M$ is isomorphic to its Sherk double cover.
\end{proposition}

%

\begin{remark}Note that $\{0\}$, $\{1\}$, $\{2\}$  generate $\cP$ under $\Delta$, and so if $\M$ is any map, then $\N = \{0\} \dub (\{1\} \dub (\{2\} \dub \M))$ has $T(\N) = \cP$, and so $\N$ is an un-double-able map in our context.\end{remark}

\begin{proposition}
\label{th:dubgp}
For any map $\M$ and any subset $I$ of $R$, $T(I\dub\M) = \left<T(\M), I\right>$
\end{proposition}
Proof:  This is obvious if $I \in T(\M)$.  If not, it is clear that $\left<T(\M), I\right>$ is contained in  $T(I\dub\M)$, and so we need only to show that  the opposite containment holds. This is equivalent to showing that if $J \in T(I\dub\M)$ and  $J \notin T(\M)$, then $I\Delta J \in T(\M)$ (see Proposition \ref{t:IDJ} ).

Let $\N  = I\dub\M$ and  note that, if $J \in T(\N)$, then $\N$ is a $J$-colorable double cover of $\M$; by Proposition \ref{th:mindub}, $\N$ must be isomorphic to $\N' = J\dub\M$.  Because $J \notin T(\M)$, there is some cycle $(f, W)$ which is $J$-inconsistent. If any such cycle is $I$-consistent, then in $\N$, the cycle $(f_0, W)$ has the same length as $(f, W)$, and so is also $J$-inconsistent, contradicting  $J \in T(\N)$.  Thus every cycle which is $J$-inconsistent is also $I$-inconsistent.  A similar argument in $\N'$ shows that every cycle in $\M$ which is $I$-inconsistent is also $J$-inconsistent. Therefore every cycle is either consistent for both, or inconsistent for both. By Lemma \ref{lemma13}, every cycle in $\M$ is $K = I\Delta J$-consistent, and then $K \in T(\M)$, as required.  \Qed

\subsection{Which maps are \texorpdfstring{$I$}{I}-double covers?}
We want to generalize the result from \cite{CW}  which says, loosely speaking,  that a reflexible map is a 2-fold orientable cover of a non-orientable reflexible map if and only if its rotation group contains an involutory element which conjugates each of the generators to its inverse.

To re-phrase that in our context, a reflexible map $\N$ which has an $R$-coloring is an $R$-double of some reflexible map $\M$ which does not have such a coloring if and only if the subgroup of  $Aut(\N)$ which preserves colors in the coloring contains an involutory element which conjugates each of the generators (of the color-preserving group) to its inverse.

If $\N$ is $I$-colorable, let $a$ be either of the $I$-colorings of $\N$ and define $C^+(\N, I)$ to be the subgroup of $C(\N)$ consisting of all $w$ such that $a(f) = a(fw)$ for all $f \in \Omega$.  Because $\N$ is $I$-colorable, this group has index 2, and so is normal, in $C(\N)$.

Note that Theorem \ref{t:isdouble} generalizes the result from \cite{CW} not only to all $I$-colorings, but also to maps with no assumption on their automorphism group.

\begin{theorem}\label{t:isdouble}
A map $\N = (\Omega, [r_0, r_1, r_2])$ is $I\dub\M$ for some non-$I$-colorable map $\M=(\Theta, [s_0, s_1, s_2])$ if and only if $\N$ is $I$-colorable and there exists $i \in I$ and  $w \in C^+(\N,I)$ of order 2   such that for every element $c \in C^+(N,I)$, $wcw=r_i c r_i$.
\end{theorem}

\textbf{Proof}: First, assume $\N= I\dub \M$ with $\M$ a non-$I$-colorable map. We will prove the existence of $w$. The group $C^+(\N,I)$ has index 2 in $C(\N)$, and there is an isomorphism $f:C^+(\N,I) \to C(\M)$ which is the letter-to-letter projection, $f(r_{i_0}r_{i_1} ... r_{i_k}) = s_{i_0}s_{i_1} ... s_{i_k}$ whenever the word $r_{i_0}r_{i_1} ... r_{i_k}$ is in $C^+(\N,I)$ and so has an even number of indices in $I$.

Because $M$ is not $I$-colorable, $I$ is not empty. So consider any $i \in I$ and let $w = f^{-1} (s_i)$. Note that $w$ is an involution since $f$ is an isomorphism. This $w$ is not $r_i$, since $r_i \notin C^+$; it is, instead, expressible as some longer product of $r_j$'s.  Then  for any $c\in C^+(\M)$ we have $f(wcw) = s_if(c)s_i$, while  $f(r_i c r_i) = s_if(c)s_i$, since $r_icr_i$ is in $C^+$, and $f$ acts there as letter-to-letter projection.  Since $f$ is an isomorphism and hence one-to-one, we conclude that $wcw = r_icr_i$.

Now assume we have $w \in C^+(\N,I)$ of order 2 and $i\in I$ such that for every element $c \in C^+(\N,I)$, $wcw=r_i c r_i$.   In particular, with $c = w$, we have $w =www = r_iwr_i$.  From this we see that $wr_i = r_iw$.  We then construct $\M$ by identifying flags in the following consistent way: let $u = w r_i$, so that $u$ is an involution.  We shall identify flag $f \in \Omega$ with $fu$. Because $u$ is an involution, this identification is unambiguous, and each flag is identified with exactly one other.  Since $w$ preserves coloring and $r_i$ changes it, we are identifying a white flag with a black flag. These equivalence classes are the flags of $\M$.

Then we have to define the connectors $s_0, s_1$ and $s_2$ of $\M$. Notice that for every $c\in C^+(\N,I)$ we have $u^{-1} c u = c$, so $u$ commutes with all of $C^+(\N,I)$. But $C^+(\N,I)$ has index 2, and $u \notin C^+(\N,I)$ so $u$ commutes with all of $C(\N)$.
Then we  define $([f, f u])s_j = [f r_j, f u r_j = f r_j u]$. The $s_i$'s are well defined and we can easily see $s_0$ and $s_2$ commute, since $r_0$ commutes with $r_2$.

To see that $\N$ is the $I$-double cover of $\M$, identify each black flag $f$ of $\N$ with the flag $[f, fu]_1$ of $I\dub\M$, and each white flag $f$ with  the flag $[f, fu]_0$. \Qed

\medskip

In order to use this theorem we would, in theory, need to prove the existence of $w$ and observe that under conjugation, $w$ acts in the same way as some $r_i$ on \emph{every} element of $C^+(\N)$. But in practice, we only need to find $w$ that acts in such a way for the \emph{generators} of $C^+(\N)$.  The element $r_i$ must not be in $C^+$ and so we will refer to it as an {\em external generator}. Table I below summarizes information about each subset $I$ of $R$. The second column (adapted from \cite{H2} ) gives generators for $C^+(\N)$ as a subgroup of $C(\N)$.  The third column gives the $r_i's$ which are not in $C^+$, and the last column shows the actions on the generators of $C^+$ under conjugation by this external generator.

\begin{table}[H]
\centering
{\small
\begin{tabular}{|c|c|c|c|} 
\hline
\multicolumn{4}{|c|}{Table I}\\
\hline
 $I$ & $C^+$ generators &external generator&action\\
\hline
$\{0\}$ & $a=r_1, b = r_2, c = r_0r_1r_0$ & $r_0$ & $a\cra c, b\cra b$\\
\hline
$\{1\}$ & $a=r_0, b = r_2, c = r_1r_0r_1, d = r_1r_2r_1$ & $r_1$ & $b \cra d, a\cra c$\\
\hline
$\{2\}$ & $a=r_0, b = r_1, c = r_2r_1r_2$ & $r_2$ & $a\cra a, b\cra c$\\
\hline
\multirow{2}{*}{$\{0,1\}$} & \multirow{2}{*}{$a=r_2, b = r_0r_1$} & $r_0$ & $a\cra a, b\cra b^{-1}$\\
\cline{3-4}
 &  & $r_1$ & $a\cra b^{-1}ab, b\cra b^{-1}$\\
\hline
\multirow{2}{*}{$\{0,2\} $}& \multirow{2}{*}{$a=r_1, b = r_0r_2, c = r_0r_1r_0$} & $r_0$ & $b\cra b, a\cra c$\\
\cline{3-4}
& & $r_2$ & $b\cra b, a\cra bcb, c\cra b a b$\\
\hline
\multirow{2}{*}{$\{1,2\}$} & \multirow{2}{*}{$a=r_0, b = r_2r_1$} & $r_2$ & $a\cra a, b\cra b^{-1}$\\
\cline{3-4}
& & $r_1$ & $a\cra b^{-1} a b, b\cra b^{-1}$\\
\hline
\multirow{3}{*}{$R$} & \multirow{3}{*}{$a=r_0r_1, b = r_1r_2$} & $r_0$ & $a\cra a^{-1}, b\cra a b^{-1} a^{-1}$\\
\cline{3-4}
& & $r_1$ & $a\cra a^{-1}, b\cra b^{-1}$\\
\cline{3-4}
& & $r_2$ & $a\cra b^{-1} a^{-1} b, b\cra b^{-1}$\\
\hline
\end{tabular}}
\end{table}

\section{Which subgroups appear on which surfaces?}\label{sec:everygroup}

The purpose of this section is to prove that every subgroup $H$ of $\cP$ occurs on every  surface on which it {\em can} appear. To clarify that statement, we call a group which includes $R$ an {\em orientable} subgroup of $\cP$, and subgroups which do not contain $R$ are {\em non-orientable}.  So we want $H$ to appear on $S$ when they are both orientable or both non-orientable. This is almost true, as explained by the following Theorem:
\begin{theorem}\label{thm:everygroup}
For every subgroup $H$ of $\cP$ and every surface $S$, there is a  map $\M$ on $S$ such that $T(\M) = H$ if and only if $S$ and $H$ are both orientable or both non-orientable, with the following exceptions:
\begin{enumerate}
\item The group $\{\emptyset, \{1\}, \{0,2\}, R\}$ does not appear on the sphere.
\item The group $\{\emptyset, \{1\}\}$  does not appear on the projective plane.
\item The group $\{\emptyset, \{0,2\}\}$ does not appear  on the projective plane.
\end{enumerate}
\label{allgp}
\end{theorem}

The proof of this claim is contained in and scattered through the remaining parts of this section.  We will consider the 16 possible groups subgroups of $\cP$ individually and in dual pairs. Such a proof would appear to require 16 constructions, but we can use Lemma \ref{DI} and Proposition \ref{th:Jup} for some simplification. Moreover, the  insertion techniques that we will introduce and the connected sum construction of Subsection \ref{sec:connectedsum} allow even more simplification. We consider the non-orientable cases first; then a simple argument, using the results of Section \ref{sec:doubles}, takes care of the orientable ones.

\subsection{\texorpdfstring{$H = \{\emptyset\}$}{\{0\}}}
Consider any triangulation of a non-orientable surface $S$. Refine the triangulation slightly by placing one new vertex inside some triangular face and connect it to the three vertices on its boundary so as to have one vertex of degree 3. Call this map $\M$. Because of the triangular faces and the the degree-3 vertex, the Cheat Sheet of Section \ref{sec:ColOri} shows us that none of the sets of size 1 or 2 are in the group, and because $S$ is non-orientable, it cannot contain $R$ either. Thus $T(\M)$ must be the trivial group.

\subsection{Insertion}\label{insert}

In this subsection, we introduce the technique of {\em insertion}.  We shall prove that any map $\M$ can be easily modified to produce a map which is face or vertex bipartite and also face or vertex pseudo-orientable.

Let us assume, for example, that we wish to modify $\M$ to construct a vertex-bipartite map. Color the vertices of $\M$ at random. For the edges for which both vertices are the same color, add a vertex (of degree 2) to split the edge in two.

The same idea can be applied if we wish vertex pseudo-orientability. We may assign orientations to the vertices at random and then split the edges whose endpoints don't match into two edges. Call this process \emph{vertex-insertion}. This can also be used to turn a map which \emph{is} vertex-bipartite or VPSO into one that is \emph{not}. Note that vertex-insertion does not modify the status of either face-bipartiteness of FPSO of $\M$.

Dually, we may replace an edge by a face bounded by two edges to get either face-bipartition or FPSO. The following picture illustrates this for face pseudo-orientation:
\begin{figure}[H]
\begin{center}
\includegraphics[height=50mm]{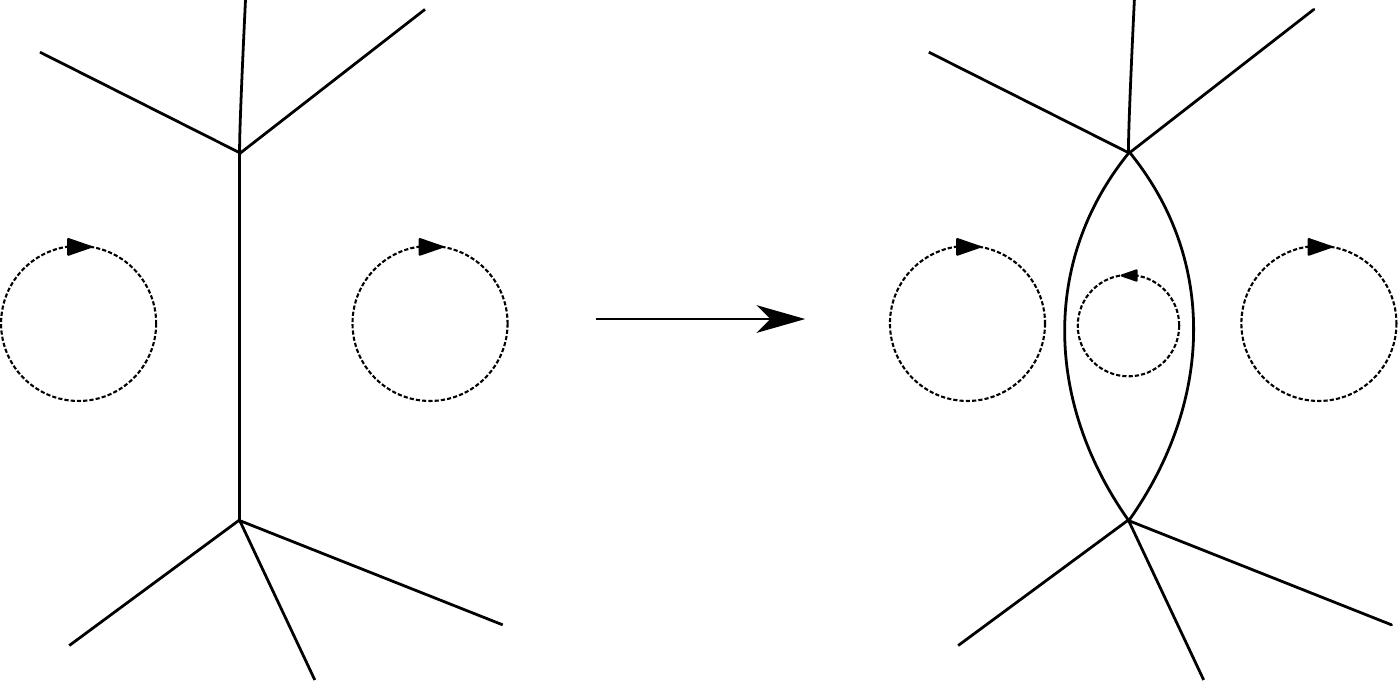}
\caption{Adjusting one edge}
\label{fig:OtoPSO}
\end{center}
\end{figure}

We begin with an arbitrary map on a surface $S$, and an arbitrary orientation of each face. If there is an edge, as in Figure \ref{fig:OtoPSO}a, where the assigned arrows flow in opposite directions along the edge, replace the edge with a face bounded by a pair of edges as shown in \ref{fig:OtoPSO}b. Orient the resulting lune with an arrow as shown.  This gives a map on the same surface with one fewer edge for which the given arrows fail to give a face pseudo-orientation. Continuing this process eventually gives us a map on $S$ admitting a face pseudo-orientation.

Let's call this process {\em face-insertion}. Here, we used it to produce a FPSO map on any surface. We can use it equally well to turn any map into a face-bipartite map on the same surface. As with vertex insertion, this same process can also be used to convert a map which {\em is} face bipartite or FPSO to one which is \emph{not}.

For many groups we can use these two operations to produce maps with that group on any non-orientable surface. We will show one of these in detail and summarize the rest.

	Consider the non-orientable group $H = \{\es, \{0,1\}\}$.   For any non-orientable surface $S$, we choose any map on $S$, and, using face-insertion as above, construct a map which is FPSO.   Now we must be certain that the group of the map is no bigger than $H$.  The other cosets of $H$ in $\cP$ are $H_2 = \{\{0\}, \{1\}\},  H_3 = \{\{0,2\}, \{1, 2\}\},  H_4 = \{\{2\}, \{0, 1, 2\}\} $.  By inserting one vertex on one edge, if necessary, we can cause the map to have a face with odd degree, eliminating $\{0\}$ and $\{1,2\}$ and hence all $H_2 \cup H_3$ from the group.
The non-orientability of the surface eliminates $R = \{0, 1, 2\}$ and hence all of $H_4$. Thus the final form $\M$ of the map has $T(\M) =  \{\es, \{0,1\}\}$, as required.

	Notice that the group $T(D(\M))$  is then   $\{\es, \{1,2 \}\}$.

We produce maps having most of the other non-orientable groups in a similar way.  Starting with an arbitrary map on a non-orientable surface we perform insertions in the following ways:
\begin{itemize}
  \item Use face-insertion to force FPSO and vertex-insertion to force VPSO. The result is a map whose group is  $\{\emptyset, \{0,1\}, \{0,2\}, \{1,2\}\}$.
  \item Use face-insertion to force FPSO and then vertex-insertion to force it to be vertex-bipartite. The resulting map has group  $\{\emptyset, \{0\}, \{1\}, \{0,1\}\}$.  Its dual has group $\{\emptyset, \{1\}, \{2\}, \{1,2\}\}$.
  \item Use face-insertion to force the map to be face-bipartite and vertex-insertion to force it to be vertex-bipartite. The resulting map has group  $\{\emptyset, \{0\}, \{2\}, \{0,2\}\}$.
  \item Use face-insertion to force face-bipartiteness and vertex-insertion in order to force some face to have an odd number of sides. The group must be just $\{\emptyset,  \{2\} \}$.  Its dual then has group $\{\emptyset,  \{0\}\}$.
\end{itemize}

\medskip

These constructions take care of all non-orientable groups except for two:  $\{\emptyset,  \{1\}\}$ and $\{\emptyset,  \{0,2\}\}$.

\subsection{The groups \texorpdfstring{$H =\{\es,\{1\}\}$}{\{0,\{1\}\}} and \texorpdfstring{$H =\{\es,\{0,2\}\}$}{\{0,\{0,2\}\}}}

In order to cope with the two non-orientable exceptional groups, we apply the technique of connected sums of maps, developed in Subsection \ref{sec:connectedsum}, as well as by giving constructions of parameterized families of suitable maps in order to have explicit constructions.

\begin{lemma}
    For each group $H \le \cP$ with $R \notin H$ and $|H| = 4$, there exists a map $\M$ on every non-orientable surface such that all the following hold:
    \begin{itemize}
        \item $T(\M) = H$.
        \item There is a face $F$ with exactly 2 sides and exactly 2 vertices.
        \item $T(\M \setminus F) = H$. That is, the flags of $\M$ that are not in $F$ are \emph{not} $I$-colorable for every $I \notin H$.
    \end{itemize}
\end{lemma}
Proof: Let $S$ be an arbitrary non-orientable surface. To construct a map on $S$ with the desired properties, first consider any map for which no face shares a vertex with itself. Because the surface is non-orientable, we may choose an arbitrarily fine subdivision for which there is a M\"obius band of faces which does not use all the faces or edges. Then apply the operations described in subsection \ref{insert} to make the map $h$-colorable for each $h$ in $H$. This is possible since any group $H$ with 4 elements with $R \notin H$ must contain two of the following: vertex-bipartite, face-bipartite, FPSO or VPSO, with the other non trivial element being the symmetric difference of the other two.

Now consider an edge that is not in any of the faces of the chosen M\"obius band and convert it into three edges, as with face-insertion but with three edges instead of two. This operation does \emph{not} change the status of vertex-bipartition, face-bipartition, VPSO or FPSO as required. Now pick either of the two newly created faces and call it $F$. Since we still have a M\"obius strip, $R \notin T(\M\setminus F)$, but every $I\in T(\M)$ is in $T(\M \setminus F)$. Therefore, $T(\M \setminus F) = H$. \Qed

\begin{proposition}
    For $H = \{\es,\{0,2\}\}$ and $H=\{\es,\{1\}\}$ we have that for every surface (except the projective plane) there exists a map $\M$ on that surface with $T(\M)=H$.
\end{proposition}
Proof: Let $S$ be a surface with non-orientable genus $k \geq 2$. Consider maps as in the previous Lemma $\M_1$ and $\M_2$ such that $\M_1$ lies on a surface with non-orientable genus $k-1$ and $\M_2$ lies on the projective plane and such that:
\begin{itemize}
    \item If $H=\{\es,\{1\}\}$, then
        $$T(\M_1) = \{\emptyset, \{1\}, \{2\}, \{1,2\}\}\text{ and }T(\M_2) = \{\emptyset, \{0\}, \{1\}, \{0,1\}\}.$$
    \item If $H=\{\es,\{0,2\}\}$, then
        $$T(\M_1) = \{\emptyset, \{0,1\}, \{0,2\}, \{1,2\}\}\text{ and }T(\M_2) = \{\emptyset, \{0\}, \{2\}, \{0,2\}\}.$$
\end{itemize}
Then $\M_1 \oplus \M_2$ (along the 2-sided faces as in the previous lemma) lives on the non-orientable surface with genus $k$ and $T(\M_1 \oplus \M_2) = H$ (see Theorem \ref{t:plus} ). \Qed

\begin{figure}[H]
\begin{center}
\includegraphics[height=60mm]{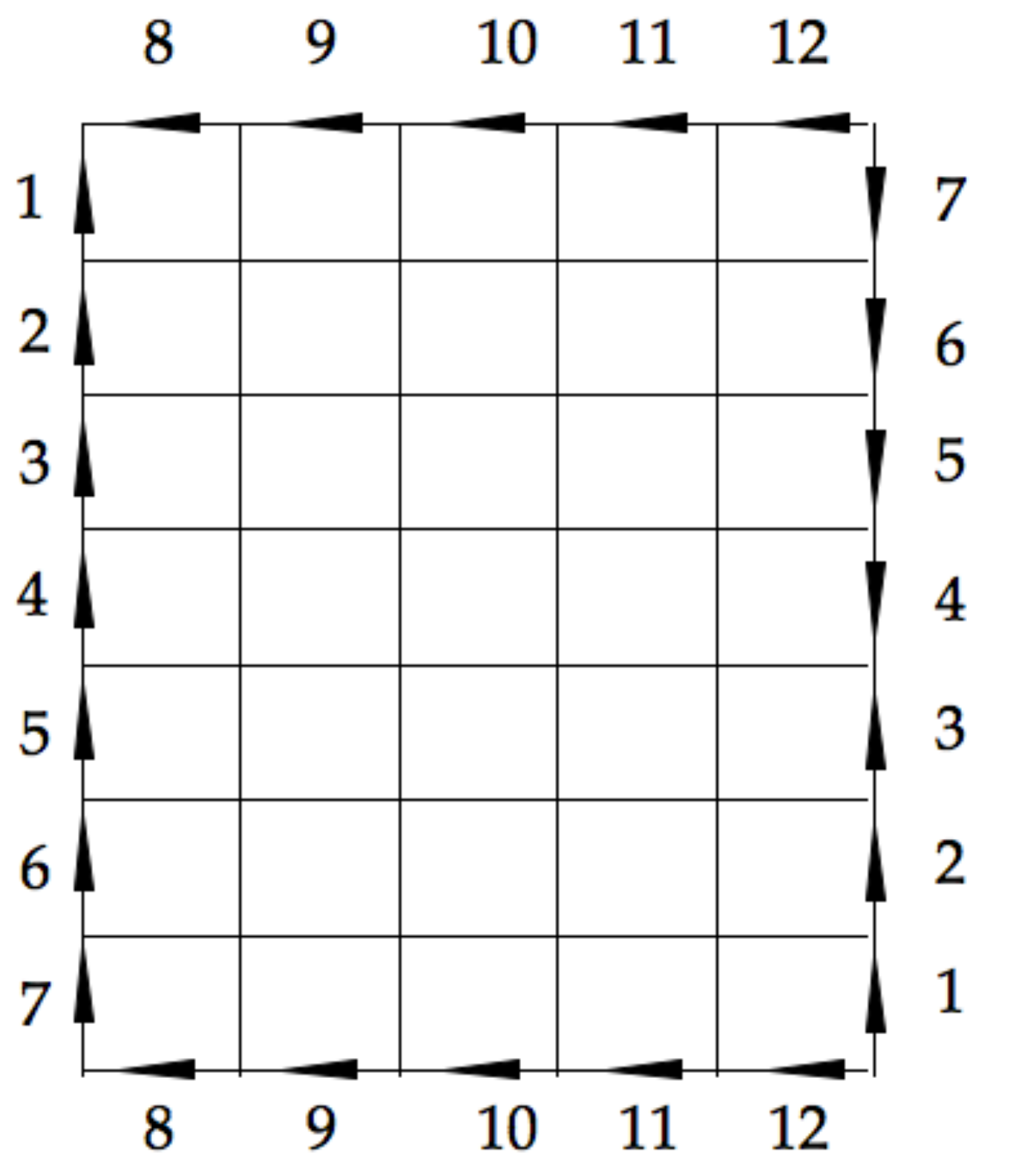}
\caption{The map $\G(5, 7, 3)$}
\label{fig:G573}
\end{center}
\end{figure}

We give explicit constructions for these last two groups:

Let $H=\{\emptyset, \{1\}\}$, and consider the map $\G = \G(m, n, k)$ where $k, m, n > 0$. This map has $mn$ faces, each one a square.   They are arranged in an $m\times n$ rectangle. The top $m$ edges are identified with the bottom $m$ directly and orientably.  Each vertical edge is identified with the one diametrically across from it; the first $k$ orientably, the rest non-orientably.   Figure \ref{fig:G573} shows $\G(5, 7, 3)$.

Because vertical and horizontal edges alternate around every face and every vertex, this map is always edge-bipartite, i.e.,$\{1\}$-colorable, and lies on the non-orientable surface of genus $k+2$. ( $\G(m, n, 0)$, then, is on the Klein bottle)    If $n$ is odd and $m  > 2$, then the map is neither face-bipartite nor vertex-bipartite, and so its group is $\{\emptyset, \{1\}\}$.

For $H=\{\emptyset, \{0,2\}\}$ Consider the maps shown in Figure \ref{fig:One}:

\begin{figure}[H]
\begin{center}
\includegraphics[height=50mm]{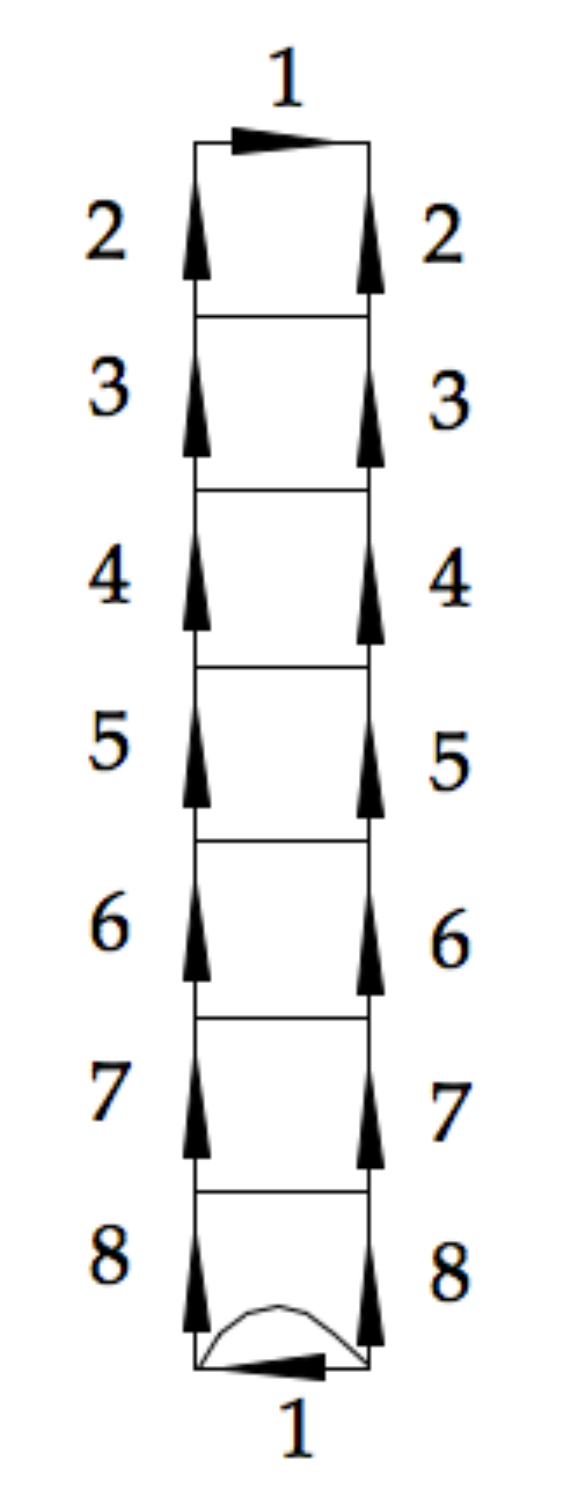}
\hspace{10mm}
\includegraphics[height=50mm]{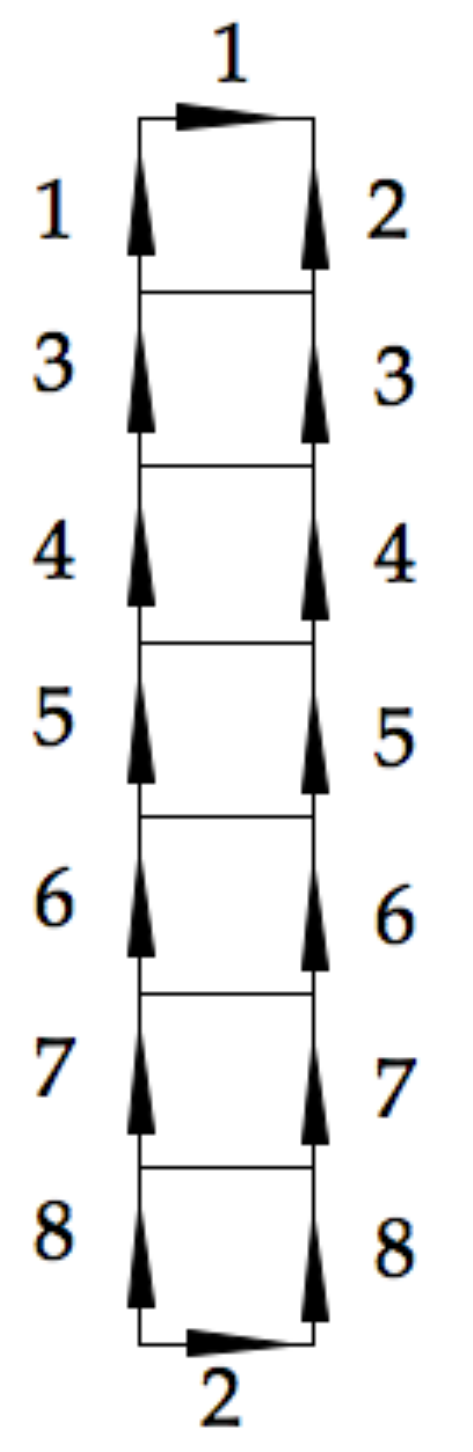}
\caption{Two maps}
\label{fig:One}
\end{center}
\end{figure}

In general we have a rectangle one unit wide, having any height. The maps shown in Figure \ref{fig:One} have characteristic 0 and 1 respectively.  In the left picture, by making the strip long enough,  we may switch enough pairs of adjacent labels on the right side (for example, switch labels 3 and 4 on the right) to make any even genus non-orientable surface. Similarly, for the picture on the right we may switch enough pairs of adjacent labels in a sufficiently long strip to construct any odd genus non-orientable surface.

To see that the group is indeed $\{\emptyset,\{0,2\}\}$, we see that neither map is face-bipartite, since some faces are glued to themselves, or face pseudo-orientable, since both contain faces that are glued to themselves in an orientable way. Similarly, neither are vertex-bipartite or vertex pseudo-orientable.

This concludes the proof for the non-orientable case.

\subsection{Orientable Groups}
If the subgroup $H$ contains $R$, then it is generated by $R$ and some subgroup $H'$ not containing $R$.  Each orientable surface $S$ is the orientable twofold cover of some non-orientable surface $S'$.  By the preceding subsections, $S'$ has a map $\M'$ whose group is $H'$ (unless S is the sphere, $S'$ is the projective plane and $H$ is $\{\emptyset,\{1\},\{0,2\},\{0,1,2\}\}$).  Then $R\dub \M'$ is a map on $S$ whose group is $H$. By Proposition \ref{th:dubgp},  $T(\M)$ is exactly $H$.  Thus, with the one exception,  for each subgroup $H$ containing $R$, every orientable surface admits a map whose group is $H$.\Qed

All that is left to prove now is that the three exceptions mentioned in the theorem are in fact exceptions -- namely that there is no map on the specified surfaces with the specified group. We shall prove this first  for the group $\{\emptyset, \{1\}, \{0,2\}, R\}$ on the sphere.

\begin{lemma}\label{l:nosphere}
  There is no map on the sphere with group $\{\emptyset, \{1\}, \{0,2\}, R\}$.
\end{lemma}
Proof: Suppose $\M$ is one such map. Since $\{0\} \notin T(\M)$, the underlying graph of $\M$ is not bipartite, which means it has an odd cycle $C$. By Jordan's Closed Curve Theorem, this divides the sphere into two connected components. Choose either and erase every edge and vertex exterior to it in order to form the deleted part into a single face whose boundary is $C$.  This gives a new map $\M'$, also on the sphere. Since $\{1\} \in T(\M)$, every face of $\M$ must have an even number of sides. Therefore $\M'$ has exactly one face with an odd number of sides. Then the underlying graph of $D(\M')$ has exactly one vertex of odd degree, a well-known impossibility. \Qed

\begin{corollary}
No map $\M$ on the projective plane has $T(\M) = \{\emptyset, \{1\}\}$ or $T(\M) = \{\emptyset, \{0,2\}\}$.
\end{corollary}
Proof: If such a map were to exist, its orientable double cover would violate Lemma \ref{l:nosphere}.\Qed

The proof of Theorem \ref{allgp} is now complete.

\begin{corollary}
No map on a non-orientable surface has an $R$-coloring.  In all other cases, for any subset $I$ of $R$, and any surface $S$, there is a map $\M$ on $S$ which is $I$-colorable.
\end{corollary}

\section{Maniplexes and higher dimensions}
\label{sec:mpx}
To generalize the results of this paper to higher-dimensional objects, we use the idea of a {\em maniplex}, first introduced in \cite{Wman}.  We summarize the definitions here:  an {\em $n$-maniplex} is a pair $(\Omega, [r_0, r_1, r_2, \dots, r_n])$, where $\Omega$ is a set of things called {\em flags}, and each $r_i$ is a fixed-point-free involution on $\Omega$ such that $\langle r_0, r_1, \dots, r_n\rangle$ is transitive on $\Omega$, and for every $i$, $j$ such that $0\le i < j-1<n$, the permutations $r_i$ and $r_j$ commute and are disjoint.

Maniplexes clearly generalize maps and slightly generalize (the flag-graphs of) abstract polytopes. In particular, every map is a 2-maniplex and every 2-maniplex is a map.

The language of maniplexes is the language of polytopes:  an $i$-face is an orbit under the group $R_i = \langle r_0, r_1, \dots, \hat{r_i}, \dots, r_n\rangle$, generated by all of the $r_j$'s except $r_i$.  A 0-face is a vertex, a 1-face is an edge, an $n$-face is a {\em facet}.

The definitions of {\em coloring, projection, cover} generalize easily. There is no notion of surface or manifold which applies to maniplexes; nevertheless, we describe a maniplex which has an $R = \{0,1,2,\dots, n\}$-coloring as being orientable. In view of Theorem \ref{thm:super}, one may define a maniplex to be $I$-pseudo-orientable if and only if it is $(R\setminus I)$-colorable. We may give a definition of FPSO that is more in the spirit of our intuition: First, we need each facet to be orientable and second, orientations can be assigned to the facets so that they agree on every subfacet. It is easy to see that this is equivalent to an $(R \setminus \{n\})$-coloring of $\M$.

The idea of an $I$-coloring, for any $I \subseteq R$, comes forward easily, as does the fact that the set of $I$'s for which a given maniplex $\M$ is $I$-colorable forms a subgroup $T(\M)$ of the power set $\cP$ of $R$ under $\Delta$.  We can define $i$-face bipartiteness: An $i$-face is a connected component under all of the $r_j$'s except $r_i$. Make these the vertices of a graph, and join two of them by an edge is some flag of one is $r_i$ adjacent to some flag of the other.  If this graph is bipartite, we say that $\M$ is $i$-face bipartite. Then, it is easy to show then that $\M$  is $i$-face bipartite if an only if it is $\{i\}$-colorable.

We can define connected sum of two $n$-maniplexes, removing isomorphic facets, one from each maniplex and adjusting the $r_n$ connections to form the sum.  The theorem that the sum has any coloring common to both still holds.

The operators generalize with a little care.  We form $opp(\M)$ from $\M$ by replacing $r_2$ with the product $r_0r_2$.  We form $D(\M)$ from $\M$ by reversing the order of the generators.  And then $P(\M) = D(opp(D(\M)))$.  Similar facts hold about how colorings of $\M$ are inherited by its direct derivates.  The medial operation is more difficult to generalize.

The construction of $I$-doubles is straightforward to generalize, and all of the facts about the doubles do as well. The fact that maps have rank 3 plays no role in any of these proofs.

Again, because there is no  notion equivalent to that of `surface', there is no natural way to generalize the results of Section \ref{sec:everygroup}.

\section{Open Questions}

The constructions in Section \ref{sec:everygroup} lead to quite general maps, although some of them may be somewhat degenerate in the viewpoint of other previous work. It may be the case that more exceptions in Theorem \ref{allgp} are needed if we ask the extra requirement that the maps are polyhedral maps (that is, the intersection of the closure of two distinct faces is either empty, a single vertex, or a single edge. See \cite{JV} ). However, this is beyond the scope of this paper.

In previous sections we considered colorings of the flags of maps insisting that we use only two colors. In doing so, we ensure that if a given flag is $i$-adjacent to another flag with the same color then all flags are the same color as their $i$-adjacent flags.

Consistent colorings with $k$ colors can be defined following the idea of \cite{B} and \cite{OPW}, that is, each flag is assigned a color in $\{1, \dots, k\}$ with the restriction that for every $i\in \{0, 1, 2\}$ and $j \in \{1, \dots, k\}$, the color of the $i$-adjacent flag to a flag $f$ colored $j$ depends on $j$ but not on $f$.

\begin{figure}[H]
\begin{center}
\includegraphics[height=40mm]{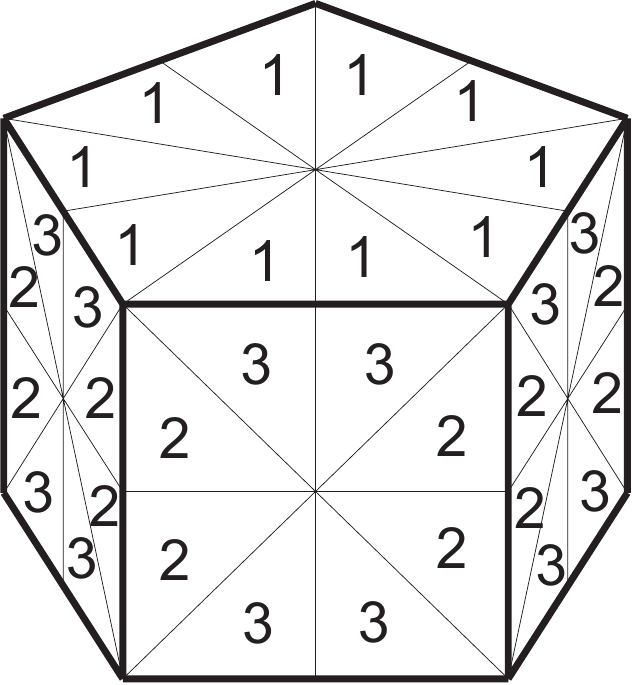}
\caption{A 3-orbit map}
\label{fig:3orb}
\end{center}
\end{figure}

When $k$-coloring a map consistently for $k \ge 3$ then it may not be true that a flag is $i$-adjacent to another flag with the same color if and only if each flag is $i$-adjacent to another flag with the same color. For example, we can color the flags of the pentagons of a pentagonal prism with color $1$, the flags of the edges shared by two squares with color $2$, and the remaining flags with color $3$, as in Figure \ref{fig:3orb}. The reader can easily verify that this is a consistent 3-coloring, however, flags colored $1$ are $1$-adjacent to flags colored $1$, and flags colored $2$ are $1$-adjacent to flags colored $3$.

It is unlikely that the concept of pseudo-orientability admits a generalization that preserves the connections between 2-colorings and pseudo-orientations of maps. However, $k$-colorings may still keep an interesting relation with $k$-fold covers of maps and maniplexes as well.

Some results in \cite{H2},  \cite{OPW}, \cite{B} relate symmetry with $I$-colorings.   We suspect that many more results can be obtained linking symmetry with $I$-colorings.

\medskip

{\bf Acknowledgements:}  The Authors, particularly the last, would like to thank the organizers of the Workshop on Abstract Polytopes (this set of organizers overlaps the set of authors) for bringing us together.  Without that meeting, this research might not have happened. The first author was partially supported by Agencija za Raziskovalno Dejavnost Republike Slovenije under the research project MU-PROM/11-010. The second and fourth authors were partially supported by PAPIIT--Mexico under grant IN112512. The third author was supported supported by DGAPA postdoctoral scholarship and additionally by PAPIIT grant IB100312.

\end{document}